\documentclass{amsart}
\usepackage{amssymb,hyperref}

\newcommand{\al}{\alpha}\newcommand{\be}{\beta}
\newcommand{\de}{\delta}
\newcommand{\ep}{\epsilon}
\newcommand{\la}{\lambda}\newcommand{\La}{\Lambda}

\newcommand{\ga}{\gamma}
\def\<{\langle}
\def\>{\rangle}
\newcommand{\R}{\mathbb{R}}\newcommand{\Z}{\mathbb{Z}}
\newcommand{\N}{\mathbb{N}}
\newcommand{\Sp}{\mathbb{S}}

\newcommand{\pt}{\partial_t}\newcommand{\pa}{\partial}

\newcommand{\D}{{\mathrm D}}

\newcommand{\beeq}{\begin{equation}}\newcommand{\eneq}{\end{equation}}

\newtheorem{thm}{Theorem}[section]
\newtheorem{prop}[thm]{Proposition}

\newtheorem{rem}{Remark}[section]

\newtheorem{lem}[thm]{Lemma}

\newenvironment{prf}{\noindent {\bf Proof.} }{\endprf\par}
\def \endprf{\hfill  {\vrule height6pt width6pt depth0pt}\medskip}

\numberwithin{equation}{section}

\begin{document}

\title[Glassey conjecture]{The Glassey conjecture with radially symmetric data}

\author{Kunio Hidano}
\address{Department of Mathematics, Faculty of Education, Mie University, 1577 Kurima-machiya-cho, Tsu, Mie 514-8507, Japan}\email{hidano@edu.mie-u.ac.jp}
\thanks{The first author was partly supported by the Grant-in-Aid for
Scientific Research (C) (No.20540165 and 23540198),
Japan Society for the Promotion of Science.}

\author{Chengbo Wang}
\address{Department of Mathematics, Zhejiang University, Hangzhou 310027, China} \email{wangcbo@gmail.com}
\thanks{The second author was supported in part by NSFC 10871175 and 10911120383.}

\author{Kazuyoshi Yokoyama}
\address{Hokkaido Institute of Technology, 7-15-4-1 Maeda, Teine-ku, Sapporo, Hokkaido 006-8585, Japan}\email{yokoyama@hit.ac.jp}


\dedicatory{} \commby{}

\begin{abstract}
In this paper, we verify the Glassey conjecture in the radial case for all spatial dimensions. Moreover, we are able to prove the existence results with low regularity assumption on the initial data and extend the solutions to the sharp lifespan. The main idea is to exploit the trace estimates and KSS type estimates.
\end{abstract}

\maketitle

\tableofcontents

\section{Introduction}

Let $n\ge 2$, $p>1$, $\Box =\pt^2-\Delta$, and $a, b$ be constants. Consider the following nonlinear wave equations
\beeq\label{eq-NLW}
    \left\{
        \begin{array}{l}
                \Box u=a |\pt u|^p+b |\nabla_x u|^p ,\ (t,x)\in \R\times\R^n\\
                u(0,x)=u_0(x)\in H^2_{\rm{rad}}(\R^n),\ \pt u(0,x)=u_1(x)\in H^1_{\rm{rad}}(\R^n)\ .
        \end{array}
    \right.
\eneq
Here $H^m_{\rm{rad}}$ stands for the space of spherically symmetric functions lying in the usual Sobolev space $H^m$.

In the 1980's, Glassey made the conjecture that the critical exponent for the problem to admit global small solutions is
$$p_c=1+\frac{2}{n-1}$$ in \cite{Glassey} (see also Schaeffer \cite{Sch86} and Rammaha \cite{Ram87}).
The conjecture has been verified
in space dimension $n=2, 3$ for general data
(Hidano and Tsutaya \cite{HiTs95} and Tzvetkov \cite{Tz98} independently)
as well as radial data (Sideris \cite{Si83} for $n=3$). For higher dimension $n\ge 4$, there are only negative results available (blow up with upper bound on expected sharp lifespan for $p\le p_c$) in Zhou \cite{Zh01}.

The purpose of this paper is to verify this conjecture in the radial case for all spatial dimensions, by proving global existence for $p>p_c$. Moreover, we are able to prove the results with low regularity assumption on the initial data and extend the solutions to the sharp lifespan (for all $1<p<1+2/(n-2)$).

Before presenting our main results, let us first give a brief review of the history.
The problem is scale invariant in the Sobolev space $\dot H^{s_c}$ with $$s_c=\frac{n}2+1-\frac{1}{p-1}\ .$$ For local well-posedness of the problem, it has been intensively studied at least for $p\in \N$, when the general result requires the initial data lie in $H^s\times H^{s-1}$ for $s>\max(s_c, (n+5)/4)$ (see Ponce and Sideris \cite{PoSi93}, Tataru \cite{Ta99}, Fang and Wang \cite{FaWa05} and references therein). If $p\ge 3$ or $p=2$ with $n\ge 4$, the problem is locally well-posed in $H^s\times H^{s-1}$ for $s>s_c$, when the initial data have radial symmetry or certain amount of angular regularity (see Fang and Wang \cite{FW10} and references therein).

For the long time existence of the solutions with $C_0^\infty$ small data of size $\ep>0$, it is well known also for the case of $p\in \N$ (even for the problem of quasilinear equations). When $p>p_c$, we have global existence. For $p=p_c$, we have almost global existence with lifespan $T_\ep$ which satisfies $$\log (T_\ep)\sim \ep^{1-p}\ .$$ Instead, if $p<p_c$, we have long time existence with lifespan $$T_\ep\sim \ep^{-\frac{p-1}{1-(n-1)(p-1)/2}},$$ see John and Klainerman \cite{JoKl84_01}, Klainerman \cite{Kl85_01}, Sogge \cite{So08} and references therein.
Moreover, the estimate on the lifespan $T_\ep$ is sharp for the problem with nonlinearity $|\pt u|^p$ (see Rammaha \cite{Ram97} for $p=2$ and $n=2,3$, Zhou \cite{Zh01} for $p\in\R$ and $1<p\le p_c$).

There is not much work on the long time existence with low-regularity small data.
In \cite{HY06}, Hidano and Yokoyama proved almost global existence for small $H^2_{\rm{rad}}\times H^1_{\rm{rad}}$ data when $p=2$ and $n=3$. It was generalized to the quasilinear problem in our recent work \cite{HWY10}.
If $p\ge 3$ or $p=2$ with $n\ge 4$, we have global (almost global for $p=3$ and $n=2$) in $H^s$ with $s>s_c$ and certain angular regularity (Sterbenz \cite{St07} and Fang and Wang \cite{FW10}).

We will use $\La_i$ to denote the norm of the initial data,
$$\La_i:=\|u_0\|_{\dot H^i(\R^n)}+\|u_1\|_{\dot H^{i-1}(\R^n)}, \ i=1,2\ .$$
Let $\pa=(\pa_x,\pt)$ with $\pa_x=(\pa_{x_1},\pa_{x_2},\dots,\pa_{x_n})$,  $x=r\omega$ with $r=|x|$ and $\omega\in\Sp^{n-1}$, and $\<r\>=\sqrt{1+r^2}$.
Now we are ready to state our main results. The first result is the global existence theorem for $p>p_c$ and $n\ge 3$.
\begin{thm}\label{thm-sub}
Let $n\ge 3$ and $1+2/(n-1)<p<1+2/(n-2)$. Consider the nonlinear wave equation \eqref{eq-NLW}. For any choice of $s_1$, $s_2$ such that
$1/2\le s_1< n/2-1/(p-1)<s_2\le  1$,
 there exist constants $C, \ep_0>0$, such that if
$$ \La_1^{1-s_1}\La_2^{s_1}+\La_1^{1-s_2}\La_2^{s_2}\le \ep_0\ ,$$
then we have a unique global solution $u$ to \eqref{eq-NLW} satisfying
$$ u\in C([0,\infty); H^2_{\rm{rad}}(\R^n))\cap C^1([0,\infty); H^1_{\rm{rad}}(\R^n))\ ,$$
 $$\|\pa u\|_{L^\infty([0,\infty); L^2(\R^n))}+\|r^{-\delta}\<r\>^{-1/2+\delta'} \pa u\|_{L^2([0,\infty)\times\R^n)}\le C \La_1\ ,$$
 $$\|\pa\pa_x u\|_{L^\infty([0,\infty); L^2(\R^n))}+\|r^{-\delta}\<r\>^{-1/2+\delta'} \pa\pa_x u\|_{L^2([0,\infty)\times\R^n)}\le C \La_2\ ,$$
where
 \beeq\label{eq-sub-delta}\de=\frac{n-2s_2}4 (p-1)\ ,\ \de'=\frac{1-(s_2-s_1)(p-1)}2\ .\eneq
\end{thm}
In contrast, when $p=p_c$, we have the almost global existence.
\begin{thm}\label{thm-crit}
Let $n\ge 3$ and $p=1+2/(n-1)$. Consider the nonlinear wave equation \eqref{eq-NLW}. For any choice of $s$ such that
$1/2< s\le  1$,
 there exist constants $C, c, \ep_0>0$, such that if
$$\ep:= \La_1^{1/2}\La_2^{1/2}+\La_1^{1-s}\La_2^{s}\le \ep_0\ ,$$
then we have a unique almost global solution $u$ to \eqref{eq-NLW} satisfying
$$ u\in C([0,T_*]; H^2_{\rm{rad}}(\R^n))\cap C^1([0,T_*]; H^1_{\rm{rad}}(\R^n))\ ,$$
 $$\|\pa u\|_{L^\infty([0,T_*]; L^2(\R^n))}+\ep^{(p-1)/2}\|r^{-\delta}\<r\>^{-1/2+\delta} \pa u\|_{L^2([0,T_*]\times\R^n)}\le C \La_1\ ,$$
 $$\|\pa\pa_x u\|_{L^\infty([0,T_*]; L^2(\R^n))}+\ep^{(p-1)/2}\|r^{-\delta}\<r\>^{-1/2+\delta} \pa\pa_x u\|_{L^2([0,T_*]\times\R^n)}\le C \La_2\ ,$$
where
$$\de=\frac{n-2s}4 (p-1),\ T_*=\exp(c \ep^{1-p})\ .$$
\end{thm}

For the case $1<p<p_c$, we expect a long time existence of the solution.
\begin{thm}\label{thm-super}
Let $n\ge 2$ and $1<p<1+2/(n-1)$. Consider the nonlinear wave equation \eqref{eq-NLW}.
There exist constants $C, c>0$, such that
we have a unique solution $u$ to \eqref{eq-NLW} satisfying
$$ u\in C([0,T_*]; H^2_{\rm{rad}}(\R^n))\cap C^1([0,T_*]; H^1_{\rm{rad}}(\R^n))\ ,$$
 $$\|\pa u\|_{L^\infty([0,T_*]; L^2(\R^n))}+T_*^{\de-1/2}\|r^{-\delta} \pa u\|_{L^2([0,T_*]\times\R^n)}\le C \La_1\ ,$$
 $$\|\pa\pa_x u\|_{L^\infty([0,T_*]; L^2(\R^n))}+T_*^{\de-1/2}\|r^{-\delta} \pa\pa_x u\|_{L^2([0,T_*]\times\R^n)}\le C \La_2\ ,$$
where
$$T_* = c (\La_1^{1/2}\La_2^{1/2})^{-\frac{2(p-1)}{2-(n-1)(p-1)}}\ ,$$
$$\de=\left\{\begin{array}{ll}
  \frac{(n-1)(p-1)}2,& 1<p<1+\frac{1}{n-1}\\
  \frac{(n-1)(p-1)}4,& 1+\frac{1}{n-1}\le p<1+\frac{2}{n-1}=p_c\ .
\end{array}\right.$$
\end{thm}
  As can be observed from the statement, for $p<p_c$, we do not require the smallness of the initial data, in contrast to $p\ge p_c$.

\begin{rem}
In our Theorems, the assumptions posed on the initial data are of ``multiplicative form", which is considered as one of the main innovations in this paper. For example, in Theorem \ref{thm-super}, the quantity $\La_1^{1/2}\La_2^{1/2}$ is in fact scale-invariant, and it scales like the homogeneous Sobolev space $\dot H^{3/2}$. The assumptions in the other two Theorems are almost critical, which scale like
$\dot H^{3/2}\cap \dot H^{3/2+\ep}$ for $p=p_c$ and $\dot H^{s_c-\ep}\cap \dot H^{s_c+\ep}$ for $p>p_c$, with the critical scaling regularity $s_c=n/2+1-1/(p-1)$. One of the advantages of using the ``multiplicative form" is that, for $p\ge p_c$, even if $\La_1$ is not so small, we still have {\rm (}almost{\rm)} global solutions when $\La_2$ is sufficiently small.
\end{rem}
Here, we would like to point out an interesting similarity between the Glassey conjecture and the Strauss conjecture. Recall that for the Strauss conjecture, where the nonlinearity is $|u|^p$, we find similar phenomena. Besides the critical regularity $s_c=n/2-2/(p-1)$, there is one more Sobolev regularity, namely $s_d=1/2-1/p$ (see Sogge \cite{So08} Section IV.4), as far as the radially symmetric functions are concerned. The critical exponent $p=p_0$ for this problem to have global small solutions is given by the positive root of the equation
$$(n-1)p^2-(n+1)p-2=0\ .$$
There is an interesting relation between these two facts: if $p>1$, we see that
$$s_c>s_d \ {\rm if\ and\ only\ if}\ p>p_0\ ,$$
and the sharp lifespan for $1<p<p_0$ has the order $\ep^{1/(s_c-s_d)}$.

Interestingly enough, for the Glassey conjecture, the index $3/2$ plays the same role as $s_d$. We have $s_c>3/2$ if and only if $p>p_c$ for $p>1$, and the sharp lifespan
$T_*$ has also the order $\ep^{1/(s_c-3/2)}$ for $p<p_c$.
 These observations strongly suggest that, for the equation \eqref{eq-NLW}, by adding certain amount of angular regularity if necessary, the minimal regularity for the problem to be well-posed is
$$\max\left(\frac{3}{2}, s_c\right)\ .$$

When $n=2$ and $p\ge p_c=3$, it seems to us that the methods to prove the preceding theorems are not sufficient to give satisfactory results. In spite of that, we can use the generalized Strichartz estimates of Smith, Sogge and Wang \cite{SSW10} to prove the following global result for $p>p_c$.
\begin{thm}\label{thm-sub-n2}
Let $n=2$ and $p>3$. Consider the nonlinear wave equation \eqref{eq-NLW}.
There exist constants $C, \ep_0>0$, such that
if $$\ep:= \Lambda_1^{1/(p-1)}\La_2^{1-1/(p-1)}\le\ep_0\ ,$$
then we have a unique global solution $u$ to \eqref{eq-NLW} satisfying
$$u\in C_t H_{\rm{rad}}^2\cap C^1_t H_{\rm{rad}}^1,\   \|\pa u\|_{L^\infty_t L^2_x}\le C \La_1,\  \|\pa \pa_x u\|_{L^\infty_t L^2_x}\le C \La_2,\   \|\pa u\|_{L^{p-1}_t L^\infty_x} \le C \ep \ .$$
\end{thm}
\begin{rem}
  For $p=p_c$ and $n=2$, it has been proved in Fang and Wang \cite{FW10} that the problem has a unique almost global solution with almost critical regularity for small data, which is not necessarily radial. A similar result for $p>3$ and $p\in \N$ has also been obtained there.
\end{rem}

This paper is organized as follows. At the end of this section, we list our basic notation. In the next section, we give several Sobolev type estimates related with the trace estimates. In Section 3, we prove some space-time $L^2$ estimates, which are variants of the Morawetz-KSS estimates.
In Sections 4 and 5, we give the proof of the (almost) global results for $n\ge 3$ (Theorems \ref{thm-sub} and \ref{thm-crit}) and the scale-supercritical result for $n\ge 2$ (Theorems \ref{thm-super}), based on the results from Sections 2 and 3. In the last section, a simple proof for $p>p_c$ and $n=2$ (Theorem \ref{thm-sub-n2}) is provided, by using the generalized Strichartz estimates of \cite{SSW10}.

{\bf Notation.} Let $\delta\in (0,1/2)$, $\de'<\de$.
We denote $\D=\sqrt{-\Delta}$ and the homogeneous Sobolev norm
$$\|u\|_{\dot H^s}=\|\D^s u\|_{L^2(\R^n)}\ .$$ The homogeneous Sobolev space $\dot H^s$ with $s<n/2$ is defined as the completion of $C_0^\infty$ with respect to the semi-norm $\|\cdot\|_{\dot H^s}$.

 For fixed $T>0$, we will use the following notation. We use $\|\cdot\|_{E_i}$ ($i=1,2$) to denote the energy norm of order $i$,
$$\|u\|_{E}=\|u\|_{E_1}=\|\pa u\|_{L^\infty([0,T];L^2(\R^n))}\ ,$$
$$\|u\|_{E_2}=\|\pa_x \pa u\|_{L^\infty([0,T];L^2(\R^n))}\ .$$
We will use $\|\cdot\|_{LE}$ to denote the local energy norm,
$$\begin{array}{ll}
\|u\|_{LE}=\|u\|_{LE_1}=&\|r^{-\delta}\<r\>^{-1/2+\delta'} \pa u\|_{L^2([0,T]\times\R^n)}\\
&+\|r^{-1-\delta}\<r\>^{-1/2+\delta'}  u\|_{L^2([0,T]\times\R^n)}\\
&+(\log(2+T))^{-1/2}\left\|r^{-\delta}\<r\>^{-1/2+\delta} \left(|\pa  u|+\frac{|u|}{r}\right)\right\|_{L^2([0,T]\times\R^n)}\\
&+T^{\delta-1/2}\left\|r^{-\delta}\left(|\pa  u|+\frac{|u|}{r}\right)\right\|_{L^2([0,T]\times\R^n)}\ .
\end{array}
$$ Here, when $n\le 2$, we will assume that there are only terms about $\pa u$.
On the basis of the space $LE$, we can define
$\|u\|_{LE_2}=\|\pa_x u\|_{LE}$ and $$LE^*=r^{-\de} \<r\>^{\de'-1/2}L^2_{t,x}+(\log(2+T))^{-1/2} r^{-\de}\<r\>^{\de-1/2}L^2_{t,x}+T^{\de-1/2}r^{-\de} L^2_{t,x}\ , $$ 
where $h\in f L^2_{t,x}$ means that $h=f g$ for some $g\in L^2_{t,x}$. 
When $T=\infty$, by $LE$ norm, we mean
$$\begin{array}{ll}\|u\|_{LE}=&\|u\|_{LE_1}=\|r^{-\delta}\<r\>^{-1/2+\delta'} \pa u\|_{L^2([0,\infty)\times\R^n)}\\
&+\|r^{-1-\delta}\<r\>^{-1/2+\delta'}  u\|_{L^2([0,\infty)\times\R^n)}\\
&+\sup_{T>0}(\log(2+T))^{-1/2}\left\|r^{-\delta}\<r\>^{-1/2+\delta} \left(|\pa  u|+\frac{|u|}{r}\right)\right\|_{L^2([0,T]\times\R^n)}\\
&+\sup_{T>0}T^{\delta-1/2}\left\|r^{-\delta}\left(|\pa  u|+\frac{|u|}{r}\right)\right\|_{L^2([0,T]\times\R^n)}\ .
\end{array}
$$

\section{Sobolev type estimates}
In this section, we give several Sobolev type estimates related with the trace estimates.

First, we state a variant of the Hardy inequality.
\begin{lem}[Hardy's inequality]\label{thm-Hardy}
Let $n\ge 2$ and $0\le s\le 1$ {\rm (}$s<1$ for $n=2${\rm)}. Then we have
\beeq\label{eq-Hardy}\|r^{-s} u\|_{L^2_x}\le C \|u\|_{L^2}^{1-s}\|\pa_r u\|_{L^2}^{s} \eneq for any $u\in H^1$.
\end{lem}
\begin{prf} We only need to give the proof for $u\in C_0^\infty$.
First, we prove \eqref{eq-Hardy} for $s\ge 1/2$. Since $s<n/2$ and $s\le 1$, we have
\begin{eqnarray*}
\left\|r^{-s} u\right\|_{L^2_x}^2
            &=& \int_{\Sp^{n-1}}\int_{0}^\infty r^{-2s}|u(r\omega)|^2 r^{n-1} dr d\omega\\
            &=&\frac{1}{n-2s}\int_{\Sp^{n-1}} \int_0^\infty  |u(r\omega)|^2 \pa_r r^{n-2s} dr d\omega\\
            &=&-\frac{1}{n-2s}\int_{\Sp^{n-1}} \int_0^\infty  \partial_r (|u(r\omega)|^2) r^{n-2s} dr d\omega\\
            &\le& \frac{2}{n-2s} \int_{\Sp^{n-1}} \int_0^\infty r^{1-2s} |u| |\pa_r u| r^{n-1} dr d\omega\\
            &\le& \frac{2}{n-2s} \|r^{1-2s} u\|_{L^2_x} \|\pa_r u\|_{L^2_x}\\
            &\le &\frac{2}{n-2s} \|(r^{-s} |u|)^{(2s-1)/s} |u|^{(1-s)/s}\|_{L^2_x} \|\pa_r u\|_{L^2_x}\\
            &\le &\frac{2}{n-2s} \|r^{-s} u\|^{(2s-1)/s}_{L^2_x} \|u\|^{(1-s)/s}_{L^2_x} \|\pa_r u\|_{L^2_x}\ ,
\end{eqnarray*} where we have applied the H\"{o}lder inequality in the last step. This gives us the required estimates with $C=(2/(n-2s))^s$ and $s\ge 1/2$.
The case $s=0$ is trivial.
For $s\in (0,1/2)$, we can use the result for $s=1/2$ to prove the estimate as follows
 \begin{eqnarray*}
\left\|r^{-s} u\right\|_{L^2_x}
            &=& \left\| \left(r^{-1/2 } |u|\right)^{2s}|u|^{1-2s}\right\|_{L^2_x} \\
            &\le & \left\| \left(r^{-1/2 } |u|\right)^{2s}\right\|_{L^{1/s}_x}\left\||u|^{1-2s}\right\|_{L^{2/(1-2s)}_x} \\
            &= & \left\| r^{-1/2 } u\right\|_{L^{2}_x}^{2s}\left\|u\right\|_{L^{2}_x}^{1-2s} \\
            &\le& \left(\frac{2}{n-1}\right)^s \left\|u\right\|_{L^{2}_x}^{1-s} \left\|\pa_r u\right\|_{L^{2}_x}^{s} \ .
\end{eqnarray*}
\end{prf}

With the help of the Hardy inequality, it will be easy to prove trace estimates.
\begin{lem}[Trace estimates]\label{thm-trace}
Let $n\ge 2$.  If $1/2\le s\le 1$ {\rm (}and $s<1$ for $n=2${\rm )}, then
\beeq\label{eq-trace}
\|r^{n/2-s} u\|_{L^\infty_r L^2_\omega}\le C \|u\|_{L^2_x}^{1-s}\|\pa_r u\|_{L^2_x}^{s}
\eneq for any $u\in H^1(\R^n)$. In particular, if $s=1/2$, we have
\beeq\label{eq-trace-edpt}
\|r^{(n-1)/2} u\|_{L^\infty_r L^2_\omega}\le C \|u\|_{L^2_x}^{1/2}\|\pa_r u\|_{L^2_x}^{1/2}\ .
\eneq
\end{lem}
\begin{prf} We only need to give the proof for $u\in C^\infty_0$. The assumptions on $s$ tell us that $n-2s>0$, $0\le 2s-1\le 1$ and $2s-1<n/2$. Then by using \eqref{eq-Hardy}, we see that
 \begin{eqnarray*}
        R^{n-2s}\|u(R\omega)\|_{L^{2}_\omega}^2
            &=&-R^{n-2s}\int_{\Sp^{n-1}}\int_R^\infty \pa_r |u(r\omega)|^2 d r d\omega\\
            &\le & 2\int_{\Sp^{n-1}}\int_0^\infty  r^{n-2s} |u| |\pa_r u| dr d\omega\\
            &=&2\int_{\Sp^{n-1}} \int_0^\infty r^{1-2s} |u| |\pa_r u|  r^{n-1} dr d\omega\\
            &\le& 2\|r^{1-2s}u\|_{L^2_x}\|\pa_r u\|_{L^2_x}\\
            &\le& C \|u\|_{L^2_x}^{2-2s}\|\pa_r u\|_{L^2_x}^{2s}\ ,
 \end{eqnarray*} with $C$ independent of $R>0$.  This completes the proof.
\end{prf}

We will also need to use the following variant of the trace estimates for the proof of Theorem \ref{thm-super} in the case of $n=2$ and $2\le p<3$.
\begin{lem}\label{thm-trace-variant}
Let $n\ge 2$.  If $s\ge 0$, then
\beeq\label{eq-trace-variant}
\|r^{s} u\|_{L^\infty_r L^2_\omega}\le \sqrt{2} \|r^{s-(n-1)/2}u\|_{L^2_x}^{1/2}\|r^{s-(n-1)/2}\pa_x u\|_{L^2_x}^{1/2}\ ,
\eneq for any $u$ such that the right hand side is finite.
\end{lem}
\begin{prf}
If $u\in C_0^\infty$, this inequality follows from a simple application of integration by parts and the Cauchy-Schwarz inequality,
  \begin{eqnarray*}
    \|r^{s} u\|_{L^2_\omega}^2
        &=& r^{2s}\int_{\Sp^{n-1}}|u(r\omega)|^2 d\omega\\
        &=& -r^{2s}\int_{\Sp^{n-1}}\int_r^\infty\pa_R|u(R\omega)|^2 d R d\omega\\
        &\le &2\int_{\Sp^{n-1}}\int_r^\infty R^{2s}|u(R\omega)||\pa_R u(R\omega)| d R d\omega\\
        &\le &2\int_{\Sp^{n-1}}\int_0^\infty R^{2s-(n-1)} |u(R\omega)| |\pa_R u(R\omega)| R^{n-1}d R d \omega\\
        &\le &2 \|r^{s-(n-1)/2} u\|_{L^2_x}\|r^{s-(n-1)/2} \pa_r u\|_{L^2_x}\ .  \end{eqnarray*} Here the condition $s\ge 0$ is used to control $r^{2s}$ by $R^{2s}$.

In general, if $u\in r^{(n-1)/2-s} L^2_x$ and $\pa_x u\in r^{(n-1)/2-s} L^2_x$, we only need to construct a $C_0^\infty$ sequence which is convergent to $u$ in the corresponding norm. Define $$u_{l,m}(x):=\psi_l(x)(\rho_m*u)(x)\ ,$$
where $\psi_l(x)=\psi(x/l)$, $\rho_m(x)=m^n \rho(mx)$, $\psi, \rho \in C^\infty_0$, $\rho\ge 0$, $\int_{\R^n}\rho(x)dx=1$ and $\psi(x)\equiv 1$ for $|x|<1$.
We recall the $n$-dimensional version of (4.2) of Lemma 4.2 in our previous paper \cite{HWY10},
\beeq\label{eq-HWY}
\int_{\R^n} \frac{\rho_m(y)}{|x-y|^\al}dy \le C |x|^{-\al},\ \al<n\ ,
\eneq where the constant $C$ is independent of $m\ge 1$.

We claim that there exists a function $m=m(l)$ such that $u_{l,m(l)}\rightarrow u$ in $r^{(n-1)/2-s} L^2_x$ as $l\rightarrow \infty$. If it is true, then
we also have $(\pa_x u)_{l,m(l)}\rightarrow \pa_x u$ in $r^{(n-1)/2-s} L^2_x$.
Notice that
$$\pa_x u_{l,m}=\frac{1}{l}(\pa_x \psi)\left(\frac{x}{l}\right)(\rho_m*u)(x)+(\pa_x u)_{l,m}(x)\ .$$
For the first term, we see that
\begin{eqnarray*}
  &&\|r^{s-(n-1)/2} \frac{1}{l}(\pa_x \psi)\left(\frac{x}{l}\right)(\rho_m*u)(x)\|_{L^2_x}
  \\
  &=& \left\|\int_{\R^n} |x|^{s-(n-1)/2} \frac{1}{l}(\pa_x \psi)\left(\frac{x}{l}\right)\rho_m(x-y)u(y)dy\right\|_{L^2_x}\\
  &\le & \frac{C}{l}\left\|\||x|^{s-(n-1)/2}\rho_m^{1/2}(x-y)u(y)\|_{L^2_y} \|\rho_m^{1/2}(x-y) \|_{L^2_y}\right\|_{{L^2_x}}\\
  &\le & \frac{C}{l}\||x|^{s-(n-1)/2}\rho_m^{1/2}(x-y)u(y)\|_{L^2_y L^2_x} \\
  &\le & \frac{C}{l}\||y|^{s-(n-1)/2}u(y)\|_{L^2_y}\rightarrow 0\ ,   \end{eqnarray*} as $l\rightarrow \infty$,
where we have used the inequality \eqref{eq-HWY} with $\al=n-1-2s$ and the fact that $s>-1/2$. This gives us the convergence of $\pa_x u_{l,m(l)}$ to $\pa_x u$.

To complete the proof, it remains to prove the claim.
Observe that $$u_{l,m}-u=\psi_l(x)\left((\rho_m*u)(x)-u(x)\right)+(\psi_l(x)-1)u(x)\ .$$
For the second term, since
$$r^{s-(n-1)/2}(\psi_l(x)-1)u(x)\rightarrow 0\ \textrm{a.e.}\ x\in\R^n$$ as $l\rightarrow \infty$,
and $$|r^{s-(n-1)/2}(\psi_l(x)-1)u(x)|^2\le C |r^{s-(n-1)/2} u(x)|^2\in L^1\ ,$$
we see that, by Lebesgue's dominated
convergence theorem, $(\psi_l(x)-1)u(x)\rightarrow 0$ in $r^{(n-1)/2-s} L^2_x$ as $l\rightarrow \infty$.

We only need to control the first term $\psi_l(x)\left((\rho_m*u)(x)-u(x)\right)$. Since $r^{s-(n-1)/2} u \in L^2_x$, for any $\ep>0$, there exists a continuous function $g$ 
such that $$\textrm{supp} \ g \subset \{x\in \R^n: R_1\le |x|\le R_2\}$$ for some $0<R_1<R_2<\infty$, and
$$\|r^{s-(n-1)/2} u-g\|_{L^2}\le \ep\ .$$ To deal with
the term $\psi_l(x)\left((\rho_m*u)(x)-u(x)\right)$, we rewrite it as follows
$$\psi_l(x)\left((\rho_m*u)(x)-u(x)\right)=\psi_l(x)\left(\rho_m*(u-G)+(\rho_m*G-G)+G-u\right)\ ,$$
where $G(x):=r^{-s+(n-1)/2}g$.
We easily see that
$$\| r^{s-(n-1)/2}\psi_l(x)(G-u)\|_{L^2}\le C\|g-r^{s-(n-1)/2} u\|_{L^2}\le C\ep\ .$$
For the term involving $\rho_m*(u-G)$, we obtain
\begin{eqnarray*}
  &&\|r^{s-(n-1)/2}  \psi_l(x)\rho_m*(u-G)(x)\|_{L^2_x}
  \\
  &\le & C \left\|\int_{\R^n} |x|^{s-(n-1)/2} \rho_m(x-y)(u-G)(y)dy\right\|_{L^2_x}\\
  &\le & C \left\|\||x|^{s-(n-1)/2}\rho_m^{1/2}(x-y)(u-G)(y)\|_{L^2_y} \|\rho_m^{1/2}(x-y) \|_{L^2_y}\right\|_{{L^2_x}}\\
  &\le & C \||x|^{s-(n-1)/2}\rho_m^{1/2}(x-y)(u-G)(y)\|_{L^2_y L^2_x} \\
  &\le & C\||y|^{s-(n-1)/2}(u-G)(y)\|_{L^2_y}\le C\ep\ ,   \end{eqnarray*}
  where we have used the inequality \eqref{eq-HWY} with $\al=n-1-2 s$ and the fact that $s>-1/2$.

Finally, we consider the term involving $(\rho_m*G-G)$. Note that $G$ is a uniformly continuous function,
\begin{eqnarray*}
  |(\rho_m*G)(x)-G(x)|&=&\left|\int_{\R^n_y} \rho_m(x-y)(G(y)-G(x))dy\right|\\
  &\le& \sup_{|y-x|<C/m, x,y\in \textrm{supp} G} |G(y)-G(x)|\rightarrow 0
\end{eqnarray*} as $m\rightarrow \infty$.
Since $\textrm{supp}\ \psi_l\subset \{x\in \R^n: |x|<C l\}$,
\begin{eqnarray*}
  &&\|r^{s-(n-1)/2} \psi_l(x)((\rho_m*G)(x)-G(x))\|_{L^2}\\
  &\le& C \|r^{s-(n-1)/2}((\rho_m*G)(x)-G(x))\|_{L^2(|x|<C l)}\\
   &\le& C l^{s+1/2} \sup_{|y-x|<C/m, x,y\in \textrm{supp} G} |G(y)-G(x)|\rightarrow 0\end{eqnarray*}
  as $m\rightarrow \infty$, for any fixed $l$. This completes the proof.
\end{prf}

As we may observe, all these estimates hold for general functions. Typically, we will apply these estimates to $\pa u$, which is not radial, even if $u$ is radial. This is the main reason for us to state all the estimates above involving the $L^2_\omega$ norm.
In this way, as we can see in the following lemma, we can easily control $\pa_x u$ and $\pa_{r} u$.
\begin{lem}\label{thm-radial}
  Let $u=u(x)$ be a radially symmetric function. Then
  \beeq\label{eq-radial}|\pa_x u|= |\pa_r u|= A_{n-1}^{-1/2}\|\pa_x u\|_{L^2_\omega} \eneq with $A_{n-1}=|\Sp^{n-1}|$.
\end{lem}
The proof is just a simple calculation. Since $u$ is radial, we see $\pa_r u$ is radial. Further,
  $$\pa_x u=\frac{x}r \pa_r u,\ |\pa_x u|=|\frac{x}r| |\pa_r u|=|\pa_r u|\ ,$$
  and $$\|\pa_r u\|_{L^2_\omega}=A_{n-1}^{1/2} |\pa_r u|\ .$$
Thus,
  $$|\pa_x u|= |\pa_r u|=A_{n-1}^{-1/2}\|\pa_r u\|_{L^2_\omega}= A_{n-1}^{-1/2}\|\pa_x u\|_{L^2_\omega} \ .$$

\section{Space-time $L^2$ estimates}
In this section, we prove the space-time $L^2$ estimates, which are variants of the Morawetz-KSS estimates.

Consider the wave equation
\beeq\label{eq-lin-wave}
    \left\{
        \begin{array}{l}
                \Box u= F,\ (t,x)\in \R\times\R^n\\
                u(0,x)=u_0(x),\ \pt u(0,x)=u_1(x)\ .
        \end{array}
    \right.
\eneq

\begin{lem}[KSS type estimates]\label{thm-KSS-HY}
Let $n\ge 1$,  $0\le \de<1/2$ and $\de'<\de$.  For any solution $u=u(t,x)$ to the wave equation \eqref{eq-lin-wave}, we have the following inequality
\beeq\label{eq-KSS-HY}
    \|u\|_{E}+\|u\|_{LE}\le C (\|\pa_x u_0\|_{L^2_x}+\|u_1\|_{L^2_x}+\|F\|_{L^1_t L^2_x})\ ,
\eneq where $C$ is independent of $T>0$ and the functions $u_0\in H^1$, $u_1\in L^2$ and $F\in L^1_t L^2_x$.
\end{lem}
This is  a standard estimate now. The estimates of this type together with the application to nonlinear wave equations originate from the work of Keel, Smith and Sogge \cite{KSS02}. The variants with $LE$ norm including the homogeneous weight $r^{-\de}$ are due to Hidano and Yokoyama \cite{HY2}.
Here, for completeness, we give a proof.

\begin{prf} To begin the proof, let us recall the classical local energy estimates of Smith-Sogge (Lemma 2.2 in \cite{SS})
  \beeq\label{eq-SmithSogge}
    \|\be(x) e^{ i t \D} f\|_{L^2 ( \R\times \R^n ) }\le C_{n,\ga,\be}\|f\|_{\dot H^\ga}\ ,
  \eneq for $\be\in C_0^\infty$ and $2 \ga\le n-1$. The inequality \eqref{eq-KSS-HY} follows from this inequality with $\ga=0$ (and $\ga=1$ for $n\ge 3$), together with the energy estimate.

First, owing to the Duhamel principle and a standard scaling argument, it is enough to prove the following six inequalities
\beeq\label{eq-pf1}
    \|r^{-\de} e^{ i t \D} f\|_{L^2 ([0,1]\times \R^n)}
    \le C\| f\|_{L^2_x}\ ,\ 0\le \de<1/2\ ,
\eneq
\beeq\label{eq-pf2}
    (\log(2+T))^{-1/2}\|r^{-\delta}\<r\>^{-1/2+\delta} e^{ i t \D} f\|_{L^2 ([0,T]\times \R^n)}
    \le C\|f\|_{L^2_x}\ ,\ \de<1/2\ ,
\eneq
\beeq\label{eq-pf3}
    \|r^{-\delta}\<r\>^{-1/2+\delta'}  e^{ i t \D} f\|_{L^2(\R\times\R^n)}
    \le C\|f\|_{L^2_x}\ ,\ \de'<\de<1/2\ ,
\eneq
\beeq\label{eq-pf1-1}
    \|r^{-1-\de} e^{ i t \D} f\|_{L^2 ([0,1]\times \R^n)}
    \le C\| f\|_{\dot H^1_x}\ ,\ 0\le \de<1/2\ ,n\ge 3\ ,
\eneq
\beeq\label{eq-pf2-1}
    (\log(2+T))^{-1/2}\|r^{-1-\delta}\<r\>^{-1/2+\delta} e^{ i t \D} f\|_{L^2 ([0,T]\times \R^n)}
    \le C\|f\|_{\dot H^1_x}\ ,\ \de<1/2\ ,n\ge 3 \,
\eneq
\beeq\label{eq-pf3-1}
    \|r^{-1-\delta}\<r\>^{-1/2+\delta'}  e^{ i t \D} f\|_{L^2(\R\times\R^n)}
    \le C\|f\|_{\dot H^1_x}\ ,\ \de'<\de<1/2\ ,n\ge 3\ .
\eneq

We begin by the proof of the first three inequalities for $r\le 1$ and \eqref{eq-pf3} for $r>1$. From \eqref{eq-SmithSogge} with $\ga=0$ and $n\ge 1$, we see that
$$\| e^{ i t \D} f\|_{L^2 (\R \times \{r\le 1\})}\le C\| f\|_{L^2}\ .$$
A standard scaling argument leads us to
\beeq\label{eq-pf4}\sup_{j\in\Z}2^{-j/2}\| e^{ i t \D} f\|_{L^2 (\R\times \{r\le 2^j\})}\le C\|f\|_{L^2}\ ,\eneq
and so for any $\de<1/2$,
\begin{eqnarray*}
  \|r^{-\de} e^{ i t \D} f\|_{L^2 (\R\times \{r\le 1\})}&\le &
   C \left(2^{(1/2-\de)j} 2^{-j/2} \| e^{ i t \D} f\|_{L^2 (\R\times \{2^{j-1}<r\le 2^{j}\})}\right)_{l^2_{j:j\le 0}}\\
&\le & C\sup_{j\le 0} 2^{-j/2} \| e^{ i t \D} f\|_{L^2 (\R\times \{2^{j-1}<r\le 2^{j}\})}\\
&\le& C\| f\|_{L^2}\ .
\end{eqnarray*}
Similarly, for any $\de'<\de<1/2$, since $r\le \<r\>$, we obtain \begin{eqnarray*}
&&  \|r^{-\de}\<r\>^{\de'-1/2} e^{ i t \D} f\|_{L^2 (\R\times \{r\ge 1\})}\\
&\le&  \|r^{\de'-\de-1/2} e^{ i t \D} f\|_{L^2 (\R\times \{r\ge 1\})}\\
  &\le  &C
    \left(2^{(\de'-\de)j} 2^{-j/2} \| e^{ i t \D} f\|_{L^2 (\R\times \{2^{j-1}\le r\le 2^{j}\})}\right)_{l^2_{j:j\ge 1}}\\
&\le & C\sup_{j\ge 1} 2^{-j/2} \| e^{ i t \D} f\|_{L^2 (\R\times \{2^{j-1}\le r\le 2^{j}\})}\\
&\le& C\| f\|_{L^2}\ ,
\end{eqnarray*} which is \eqref{eq-pf3} for $r>1$.

It remains to prove \eqref{eq-pf1} and \eqref{eq-pf2} for $r>1$. For \eqref{eq-pf1}, because of the assumption $\de\ge 0$, we can easily get by the energy estimates
\begin{eqnarray*}
      \|r^{-\de} e^{ i t \D} f\|_{L^2 ([0,1]\times \{r>1\})}&\le&
            \| e^{ i t \D} f\|_{L^2 ([0,1]\times \R^n)}\\
        &\le&
            \| e^{ i t \D} f\|_{L^\infty ([0,1]; L^2(\R^n))}\\
        &\le& C\| f\|_{L^2_x}\ .
\end{eqnarray*}

For \eqref{eq-pf2} with $r>1$, we consider $1\le r\le T$ and $r\ge T$ separately.
For $r\ge T$, since $\de-1/2<0$ and $r\le \<r\>$, we obtain \begin{eqnarray*}
  \|r^{-\de}\<r\>^{\de-1/2} e^{ i t \D} f\|_{L^2 ([0,T]\times \{r\ge T\})}&\le&
  \|r^{-1/2} e^{ i t \D} f\|_{L^2 ([0,T]\times \{r\ge T\})}\\
  &\le&  T^{-1/2}\| e^{ i t \D} f\|_{L^2 ([0,T]\times \R^n)}\\
  &\le  &\| e^{ i t \D} f\|_{L^\infty ([0,T]; L^2(\R^n))}\\
        &\le& C\| f\|_{L^2_x}\ .
  \end{eqnarray*}

Now we give the estimate of \eqref{eq-pf2} for $1\le r\le T$. By \eqref{eq-pf4} and the elementary inequality $2^{[10\log(2+T)]}\ge T$ (where $[M]$ denotes the greatest integer not greater than $M$), we have
\begin{eqnarray*}
&&  \|r^{-\de}\<r\>^{\de-1/2} e^{ i t \D} f\|_{L^2 (\R\times \{1\le r\le T\})}\\
&\le&
  C \|r^{-1/2} e^{ i t \D} f\|_{L^2 (\R\times \{1\le r\le T\})}\\
  &\le&C  \left(2^{-j/2}\| e^{ i t \D} f\|_{L^2 (\R\times \{2^{j-1}\le r\le 2^{j}\})}\right)_{l^2_{j:1\le j\le 10\log (2+T)}}\\
        &\le&C (\log(2+T))^{1/2}\sup_{j} 2^{-j/2}\| e^{ i t \D} f\|_{L^2 (\R\times \{2^{j-1}\le r\le 2^{j}\})}\\
        &\le& C (\log(2+T))^{1/2} \| f\|_{L^2_x}\ .
  \end{eqnarray*} This completes the proof of the first three inequalities.

The inequalities \eqref{eq-pf1-1}-\eqref{eq-pf3-1} follow from basically the same proof, by using \eqref{eq-SmithSogge} with $\ga=1$ and Hardy's inequality. For example,
to prove \eqref{eq-pf1-1} for $n\ge 3$, we use \eqref{eq-SmithSogge} with $\ga=1$ (since $1\le (n-1)/2$ for $n\ge 3$), which tells us that
$$\| e^{ i t \D} f\|_{L^2 (\R \times \{r\le 1\})}\le C\| f\|_{\dot H^1}\ .$$
A standard scaling argument leads us to
\beeq\label{eq-pf4-1}\sup_{j\in\Z}2^{-3 j/2}\| e^{ i t \D} f\|_{L^2 (\R\times \{r\le 2^j\})}\le C\|f\|_{\dot H^1}\ ,\eneq
and so for any $\de<1/2$,
\begin{eqnarray*}
  \|r^{-\de-1} e^{ i t \D} f\|_{L^2 (\R\times \{r\le 1\})}&\le &
   C \left(2^{(1/2-\de)j} 2^{-3j/2} \| e^{ i t \D} f\|_{L^2 (\R\times \{2^{j-1}<r\le 2^{j}\})}\right)_{l^2_{j:j\le 0}}\\
&\le & C\sup_{j\le 0} 2^{-3j/2} \| e^{ i t \D} f\|_{L^2 (\R\times \{2^{j-1}<r\le 2^{j}\})}\\
&\le& C\| f\|_{\dot H^1}\ .
\end{eqnarray*}
For \eqref{eq-pf1-1} with $r>1$, since $\de\ge 0$, we can easily get by the energy estimates and Hardy's inequality \eqref{eq-Hardy} with $s=1$,
\begin{eqnarray*}
      \|r^{-1-\de} e^{ i t \D} f\|_{L^2 ([0,1]\times \{r>1\})}&\le&
            \|r^{-1} e^{ i t \D} f\|_{L^2 ([0,1]\times \R^n)}\\
        &\le&
            \|\D e^{ i t \D} f\|_{L^\infty ([0,1]; L^2(\R^n))}\\
        &\le& C\| f\|_{\dot H^1_x}\ ,
\end{eqnarray*}
which completes the proof of \eqref{eq-pf1-1}.
\end{prf}

When $n\ge 3$, we can prove the following inhomogeneous KSS type estimates with $LE^*$ norm on $F$.
\begin{lem}[
Inhomogeneous KSS type estimates]\label{thm-KSS-HY-inh}
Let $n\ge 3$, $0 < \de < 1/2$, and $\de'<\de$. For any solution $u=u(t,x)$ to the wave equation \eqref{eq-lin-wave}, we have the following inequality
\beeq\label{eq-KSS-HY-inh}
    \|u\|_{E}+\|u\|_{LE}\le C (\|\pa_x u_0\|_{L^2_x}+\|u_1\|_{L^2_x}+\|F\|_{LE^*})\ ,
\eneq where $C$ is independent of $T>0$ and the functions $u_0\in \dot H^1$, $u_1\in L^2$ and $F\in LE^*$.
\end{lem}

\begin{prf}
  i) Let us first consider smooth solutions. For such a case, we have the space-time $L^2$ estimates even for certain small perturbations of the Minkowski metric (see \cite{St05}, \cite{MeSo06_01} and our previous work \cite{HWY10}). Recall that using Lemma 2.3 and (2.30) of \cite{HWY10}, we can get
\begin{eqnarray}
 \lefteqn{T^{2\de-1}\!\!\!\int_0^T\!\!\! \int_{\{x \in {\mathbb R}^n ;\, 1 < r < T\}}\! \Bigl(\frac{|\partial u|^2}{r^{2\de}} + \frac{u^2}{r^{2+2\de}}\Bigr)dxdt}
 \label{eb24}\\
 \lefteqn{+ \bigl(\log (2+T)\bigr)^{-1}\!\!\!\int_0^T\!\!\!
\int_{\{x \in {\mathbb R}^n ;\, 1 < r < T\}}\! \Bigl(\frac{|\partial u|^2}{r} + \frac{u^2}{r^3}\Bigr)dxdt}
 \nonumber \\
 \lefteqn{+ \!\!\!\int_0^T\!\!\!
\int_{\{x \in {\mathbb R}^n ;\, 1 < r < \infty\}}\! \Bigl(\frac{|\partial u|^2}{r^{1+2\de-2\de'}} + \frac{u^2}{r^{3+2\de-2\de'}}\Bigr)dxdt}
 \nonumber \\
 &\leq& C(\|\nabla u_0\|^2_{L^2({\mathbb R}^n)} + \|u_1\|^2_{L^2({\mathbb
  R}^n)}) \nonumber \\
 & & + C\int_0^T \int_{{\mathbb R}^n} \left(
|\partial u||F| + \frac{|u||F|}{\langle r\rangle}  \right)dxdt\ ,
 \nonumber
\end{eqnarray}
for any smooth solution $u$ to the wave equations \eqref{eq-lin-wave}, $T > 1$, $\de'<\de$ and $0 < \de < 1/2$.
We will also need a slight variant of Lemma 2.2 of \cite{HWY10}. Observe that if we choose the function $$f(r)=\left(\frac{r}{r+\la}\right)^{k}$$ with $k=1-2\de\in (0,1)$ and $\la>0$, then the same argument as in the proof of Lemma 2.2 of \cite{HWY10} will tell us that
\begin{eqnarray}
\la^{2\de-1} \lefteqn{\int_0^T\!\! \int_{\{x \in {\mathbb R}^n ;\, r < \la\}} \Bigl(\frac{|\partial u|^2}{r^{2\de}} + \frac{u^2}{r^{2+2\de}}\Bigr)dxdt}
 \label{eb6}\\
 &\leq& C(\|\nabla u_0\|^2_{L^2({\mathbb R}^n)} + \|u_1\|^2_{L^2({\mathbb
  R}^n)}) \nonumber \\
 & & + C\int_0^T\!\! \int_{{\mathbb R}^n} \left(
|\partial u||F| + \frac{|u||F|}{r^{2\de}( r+\la)^{1-2\de}}
\right)dxdt,
 \nonumber
\end{eqnarray}
where the constant $C$ is independent of $\la>0$.
We only need to check the new relations (instead of (2.15) and (2.16) there)
$$\frac{f}{r}-f'(r)\ge (1-k) \frac{r^{k-1}}{(\la+r)^k},\ \Delta \left(\frac{f}{r}\right)\le -\frac{k(1-k)\la^2}{r^{3-k}(\la+r)^{2+k}} ,$$
and substitute these new relations to (2.12) and (2.17) there.

On the basis of \eqref{eb24} and \eqref{eb6}, together with the standard energy estimate
\beeq\label{eb-en}
\sup_{t\in [0,T]}\int_{\R^n}|\pa u|^2 dx \le C\left(\|\pa u(0)\|_{L^2(\R^n)}^2+\int_0^T \int_{\R^n}|\pa u| |F| dx dt\right)\ ,
\eneq it will be easy to prove the required estimates for the smooth solutions.
Suppose $T>1$ first. By applying \eqref{eb-en} to the integrals over $\{r>T\}$,
we see that
\begin{eqnarray*}
&&\|u\|_{E}^2+\|u\|_{LE(r>T)}^2\\
&\le& \sup_{t\in [0,T]}\ \!\!\! \int_{\R^n}\! |\partial u|^2 dx
+T^{2\de-1}\!\!\!\int_0^T\!\!\! \int_{\{x \in {\mathbb R}^n ;\, r> T\}}\! \Bigl(\frac{|\partial u|^2}{r^{2\de}} + \frac{u^2}{r^{2+2\de}}\Bigr)dxdt
\\
&&+ \bigl(\log (2+T)\bigr)^{-1}\!\!\!\int_0^T\!\!\!
\int_{\{x \in {\mathbb R}^n ;\,  r > T\}}\! \Bigl(\frac{|\partial u|^2}{r} + \frac{u^2}{r^3}\Bigr)dxdt
 \nonumber \\
&&+ \!\!\!\int_0^T\!\!\!
\int_{\{x \in {\mathbb R}^n ;\, r>T\}}\! \Bigl(\frac{|\partial u|^2}{r^{1+2\de-2\de'}} + \frac{u^2}{r^{3+2\de-2\de'}}\Bigr)dxdt
 \nonumber \\
&\le& \sup_{t\in [0,T]}\ \!\!\! \int_{\R^n}\! |\partial u|^2 dx+3 T^{-1}\!\!\!\int_0^T\!\!\! \int_{\{x \in {\mathbb R}^n ;\, r> T\}}\! \Bigl(|\partial u|^2 + \frac{u^2}{r^{2}}\Bigr)dxdt
\nonumber\\
&\le& C \sup_{t\in [0,T]}\ \!\!\! \int_{\R^n}\! \Bigl(|\partial u|^2 + \frac{u^2}{r^{2}}\Bigr)dx\nonumber\\
&\le& C \sup_{t\in [0,T]}\ \!\!\! \int_{\R^n}\! |\partial u|^2 dx\nonumber\\
 &\leq& C(\|\nabla u_0\|^2_{L^2({\mathbb R}^n)} + \|u_1\|^2_{L^2({\mathbb
  R}^n)} + \int_0^T \int_{{\mathbb R}^n}
|\partial u||F| dxdt ),
 \nonumber
\end{eqnarray*}
where we have applied the Hardy inequality \eqref{eq-Hardy} with $s=1$.
For the integral over $\{r<T\}$, we use \eqref{eb24} and \eqref{eb6} with $\la=1$ to get
$$
\|u\|_{LE(r<T)}^2
\leq C(\|\nabla u_0\|^2_{L^2({\mathbb R}^n)} + \|u_1\|^2_{L^2({\mathbb
  R}^n)}) + C\int_0^T \int_{{\mathbb R}^n} \left(
|\partial u||F| + \frac{|u||F|}{r}  \right)dxdt.
$$ Then an application of the Cauchy-Schwarz inequality yields the required estimate \eqref{eq-KSS-HY-inh} for $T>1$.

To prove the general result for any $T>0$, we only need to control the term
$$A[u]=T^{\delta-1/2}\left\|r^{-\delta}\left(|\pa  u|+\frac{|u|}{r}\right)\right\|_{L^2([0,T]\times\R^n)}$$ for $T\in (0,1]$.
To control this, we only need to apply  \eqref{eb6} with $\la=T$ and \eqref{eb-en} as follows:
\begin{eqnarray*}
  A[u]^2 +\|u\|_E^2& = & T^{2\delta - 1}\left\|r^{-\delta}\left(|\pa  u|+\frac{|u|}{r}\right)\right\|_{L^2([0,T]\times\R^n)}^2+\|u\|_E^2\\
   & \le & C T^{2 \delta - 1} \left\|r^{-\delta}\left(|\pa  u|+\frac{|u|}{r}\right)\right\|_{L^2([0,T]\times\R^n: r<T)}^2\\
   &  & +C T^{ - 1} \left\||\pa  u|+\frac{|u|}{r}\right\|_{L^2([0,T]\times\R^n: r>T)}^2+\|u\|_E^2\\
   & \le & C T^{2 \delta - 1} \left\|r^{-\delta}\left(|\pa  u|+\frac{|u|}{r}\right)\right\|_{L^2([0,T]\times\R^n: r<T)}^2\\
  &&+    C\left\||\pa  u|+\frac{|u|}{r}\right\|_{L^\infty_t([0,T]; L^2(\R^n))}^2\\
   &\le  & C (\|\nabla u_0\|^2_{L^2({\mathbb R}^n)} + \|u_1\|^2_{L^2({\mathbb
  R}^n)}) \nonumber \\
 & & + C\int_0^T\!\! \int_{{\mathbb R}^n} \left(
|\partial u||F| + \frac{|u||F|}{r}
\right)dxdt\ .
\end{eqnarray*}
Once again, an application of the Cauchy-Schwarz inequality gives us the required estimate \eqref{eq-KSS-HY-inh} for $T\le 1$.

ii) We next consider the case where $u$ is not smooth. By Lemma \ref{thm-KSS-HY}, we only need to prove for the case $u_0=u_1=0$. Fix $T\in (0,\infty)$. Observe that for $0<\de<1/2$, we have the Hardy inequality
$$\|r^{-\de} u\|_{L^2}\le C \|u\|_{\dot H^\delta}\le C\|u\|_{H^\delta},$$
which means $r^{-\de}L^2_x\subset H^{-\de}$, and so $LE^*\subset L^1_t H^{-\de}([0,T]\times \R^n)$ if $T<\infty$. Thus by the standard existence and uniqueness result of the linear wave equation, we have $u\in C_t H_x^{1-\de}\cap C^1_t H^{-\de}_x ([0,T]\times \R^n)$.

We claim that there exists a sequence of smooth functions $F_k$ such that $F_k\rightarrow F$ in $LE^*$. If it is true, then $u_k$ are Cauchy sequence in $E_1\cap LE_1$, and $u_k\rightarrow u$ in $C_t H_x^{1-\de}\cap C^1_t H^{-\de}_x ([0,T]\times \R^n)$. This tells us that $\{u_k\}_{k=1}^\infty$ converges to $u$ in $E_1\cap LE_1$, and so
$$\|u\|_{E_1\cap LE_1}= \lim_{k\rightarrow\infty} \|u_k\|_{E_1\cap LE_1}\le C \lim_{k\rightarrow\infty} \|F_k\|_{LE^*}=C\|F\|_{LE^*}\ ,$$
which implies \eqref{eq-KSS-HY-inh}.

To complete the proof, it remains to prove the claim.\\
{\bf Proof of the claim.}
Without loss of generality, we give the proof for $F\in r^{-\de} L^2_{t,x}$. Let $\tilde{F}(t,x)$ be the
 zero extension of $r^\de F\in L^2_{t,x}([0,T]\times\R^n)$ in $\R\times\R^n$. Let $\phi(x) \in C_0^\infty(\R^n)$ be a function with the properties $\phi\ge 0$, $\int_{\R^n} \phi(x) dx=1$, $\phi=1$ near $0$. We will also choose its one-dimensional counterpart $\psi(t)\in C_0^\infty(\R)$. Define $\phi_k(x)=2^{kn}\phi(2^k x)$ and $\psi_k(t)=2^k \psi(2^k t)$. Then the standard results of approximations of the identity give us
 $$\tilde{F}_k=(\phi_k(x)\psi_k(t))*_{t,x} \tilde F\rightarrow \tilde F\ { \rm in }\ L^2_{t,x} ([0,T]\times\R^n)\ ,$$ that is,
  $${F}^1_k := r^{-\de} \tilde{F}_k \rightarrow  F\ { \rm in }\  r^{-\de} L^2_{t,x} ([0,T]\times\R^n)\ .$$
  Notice that $F^1_k$ is smooth except at $x=0$. It suffices to set $F_k (t,x)=(1-\phi(2^k x)) F^1_k(t,x)$, which is smooth for any $t,x$.
Indeed, by Lebesgue's dominated convergence theorem, we see
\begin{eqnarray*}
 && \|r^\de (F_k-F)\|_{L^2_{t,x}([0,T]\times\R^n) } \\& =& \|\tilde F_k (1-\phi(2^k x))-\tilde F\|_{L^2_{t,x}([0,T]\times\R^n) }\\
  & \le & \|(\tilde F_k-\tilde F) (1-\phi(2^k x))\|_{L^2_{t,x}([0,T]\times\R^n) }+\|\phi(2^k x) \tilde F\|_{L^2_{t,x}([0,T]\times\R^n) }  \\
    & \le & \|\tilde F_k-\tilde F\|_{L^2_{t,x}([0,T]\times\R^n) }+\|\phi(2^k x) \tilde F\|_{L^2_{t,x}([0,T]\times\R^n) }  \rightarrow 0 \ {\rm as}\ k \rightarrow\infty\ .
\end{eqnarray*}
This completes the proof of the claim and hence that of Lemma \ref{thm-KSS-HY-inh}.
\end{prf}


\section{Glassey conjecture when $n\ge 3$}\label{sec-highD-Glassey}
Now we are ready to present our proof of the Glassey conjecture for radial initial data.

Let us first formulate the setup of the proof for existence and uniqueness.
Define
\begin{eqnarray}
  X_T:&= \{&u\in C([0,T]; H^1_{\rm{rad}}(\R^n))\cap L^\infty([0,T]; H^2_{\rm{rad}}(\R^n)):\label{eq-super-space}\\
       &&   \pt u\in C([0,T]; L^2_{\rm{rad}}(\R^n))\cap L^\infty([0,T]; H^1_{\rm{rad}}(\R^n)),\nonumber\\
       &&   \|u\|_{E_1\cap E_2}+\|u\|_{LE_1\cap LE_2}<\infty\}\ .\nonumber
\end{eqnarray}
For $R_1>0$ and $R_2>0$, we next define
$$X(R_1, R_2; T):=\{u\in X_T: \|u\|_{E_i}+\|u\|_{LE_i}\le R_i, i=1,2\}\ .$$
Endowed with
\beeq\label{eq-super-metric}
\rho(u,v):=\|u-v\|_{E_1}+\|u-v\|_{LE_1}\ ,
\eneq it is easy to check that $X(R_1, R_2; T)$ is complete with the metric $\rho(u,v)$.

For fixed $(u_0,u_1)\in H^2_{\rm{rad}}\times H^1_{\rm{rad}}$, we define the iteration map
\beeq\label{eq-Picard}
\Phi[u](t):=u^{(0)}(t)+I[N[u]]\ ,
\eneq where $u^{(0)}(t)=\cos(t\D) u_0+\D^{-1}\sin(t\D)u_1$ is the solution of the linear Cauchy problem,
\beeq\label{eq-Nonlinearity}
N[u]:=a|\pa_t u|^p+b|\nabla_x u|^p\ ,
\eneq
and
\beeq\label{eq-inhomo}
I[F]:=\int_0^t \frac{\sin((t-s)\D)}{\D} F(s) d s\ .
\eneq
For the nonlinearity $N[u]$, we have the properties \beeq\label{eq-Nonlinearity-rel1}|\pa_x^\al N[u]|\le C |\pa u|^{p-1}|\pa_x^\al \pa u|\ ,\ |\al|\le 1 ,\eneq and
  \beeq\label{eq-Nonlinearity-rel2}|N[u]-N[v]|\le C(|\pa u|^{p-1}+|\pa v|^{p-1}) |\pa (u- v)|\ .\eneq

Notice that $v=\Phi[u]$ is defined as the solution to the following equation \beeq\label{eq-lin-wave-Picard}
    \left\{
        \begin{array}{l}
                \Box v= N[u],\ (t,x)\in \R\times\R^n\\
                v(0,x)=u_0(x),\ \pt v(0,x)=u_1(x)\ .
        \end{array}
    \right.
\eneq
We aim at showing that $\Phi$ is a contraction mapping of $X(R_1, R_2; T)$, if we choose $R_1$, $R_2$ and $T$ suitably.

\subsection{Glassey conjecture when $p>p_c$ and $n\ge 3$}

Consider the nonlinear wave equation \eqref{eq-NLW} for $p>p_c$ and $n\ge 3$.

Let us begin with the estimate of the homogeneous solution, $u^{(0)}$, which follows directly from the application of Lemma \ref{thm-KSS-HY} to $u$ and $\pa_x u$ with $F=0$.
\begin{prop}\label{thm-sub-homo}
Let $n\ge 2$.  There is a positive constant $C_1$, independent of $T>0$, such that the following estimates hold
  \beeq\label{eq-sub-homo1}
   \|u^{(0)}\|_{E_1}+\|u^{(0)}\|_{LE_1}\le C_1(\|\pa_x u_0\|_{L^2}+\|u_1\|_{L^2})\ ,
  \eneq
  \beeq\label{eq-sub-homo2}
   \|u^{(0)}\|_{E_2}+\|u^{(0)}\|_{LE_2}\le C_1(\|\pa_x^2 u_0\|_{L^2}+\|\pa_x u_1\|_{L^2})\ .
  \eneq
\end{prop}

Next, we give the estimate for the inhomogeneous part.
\begin{prop}\label{thm-sub-inho} Let $p_c<p<1+2/(n-2)$ and $n\ge 3$, $u\in X_\infty$ and $s_1$, $s_2$ such that $1/2\le s_1< n/2-1/(p-1)<s_2\le  1$. Set $\de$ and $\de'$ as in \eqref{eq-sub-delta}.
 Then
  there is a positive constant $C_2$, such that the following estimates hold
\beeq\|I[N[u]]\|_{E_i\cap LE_i}\le C_2 (\|u\|_{E_1}^{1-s_1}\|u\|_{E_2}^{s_1}+\|u\|_{E_1}^{1-s_2}\|u\|_{E_2}^{s_2})^{p-1} \|u\|_{LE_i},\ i=1,2\ .\label{eq-sub-inhomo}
\eneq
  Moreover, if $u,v\in X_\infty$, we have
\begin{eqnarray}
&&\|\Phi[u]-\Phi[v]\|_{E_1\cap LE_1} \label{eq-sub-converg}\\
&  &\le C_3 \|u-v\|_{LE_1} \nonumber\\
&&\times (\|u\|_{E_1}^{1-s_1}\|u\|_{E_2}^{s_1}+\|u\|_{E_1}^{1-s_2}\|u\|_{E_2}^{s_2} +\|v\|_{E_1}^{1-s_1}\|v\|_{E_2}^{s_1}+\|v\|_{E_1}^{1-s_2}\|v\|_{E_2}^{s_2})^{p-1}\ ,  \nonumber
\end{eqnarray}
 for some $C_3$.
\end{prop}

\begin{prf}
First, by Lemma \ref{thm-trace}, we have for any $s\in [1/2,1]$,
$$\|r^{n/2-s} u\|_{L^\infty_r L^2_\omega}\le C \|u\|_{L^2_x}^{1-s}\|\pa_r u\|_{L^2_x}^{s}\ . $$ Fix $s_1$, $s_2$ such that $1/2\le s_1< n/2-1/(p-1)<s_2\le 1$. Then we have
\beeq\label{eq-sub-decay} \|u(r\omega)\|_{L^2_\omega}\le C r^{s_2-n/2}\<r\>^{s_1-s_2}(\|u\|_{L^2_x}^{1-s_1}\|\pa_r u\|_{L^2_x}^{s_1}+\|u\|_{L^2_x}^{1-s_2}\|\pa_r u\|_{L^2_x}^{s_2})\ .\eneq
By Lemma \ref{thm-radial}, we have for $u\in X_\infty$,
  \beeq\label{eq-sub-decay2}
    |\pa u|\le C r^{s_2-n/2}\<r\>^{s_1-s_2}(\|u\|_{E_1}^{1-s_1}\|u\|_{E_2}^{s_1}+\|u\|_{E_1}^{1-s_2}\|u\|_{E_2}^{s_2})\ .
  \eneq
From \eqref{eq-sub-decay2}, \eqref{eq-Nonlinearity-rel1} and \eqref{eq-sub-delta}, it is clear that, for $i=1,2$,
  \begin{eqnarray*}
 &&   \sum_{|\al|= i-1}\|r^\de\<r\>^{1/2-\de'}\pa_x^\al N[u]\|_{L^2_x}\\
  &\le& C \sum_{|\al|= i-1}\|r^\de\<r\>^{1/2-\de'}|\pa u|^{p-1}\pa_x^\al \pa u\|_{L^2_x} \\
          &\le & C  (\|u\|_{E_1}^{1-s_1}\|u\|_{E_2}^{s_1}+\|u\|_{E_1}^{1-s_2}\|u\|_{E_2}^{s_2})^{p-1}\\
&&          \times\sum_{|\al|= i-1}\|r^{\de+(s_2-n/2)(p-1)}\<r\>^{ 1/2-\de'+(s_1-s_2)(p-1)}\pa_x^\al \pa u\|_{L^2_x}\\
          &= & C (\|u\|_{E_1}^{1-s_1}\|u\|_{E_2}^{s_1}+\|u\|_{E_1}^{1-s_2}\|u\|_{E_2}^{s_2})^{p-1}\sum_{|\al|= i-1}\|r^{-\de}\<r\>^{-1/2+\de'}\pa_x^\al \pa u\|_{L^2_x}\ .\\
  \end{eqnarray*}
It is easy to check that $\de$ and $\de'$ satisfy $0<\de<1/2$ and $\de'<\de$. Now applying Lemma \ref{thm-KSS-HY-inh} to $\pa_x^\al u$ with $|\al|\le 1$ and $u_0=u_1=0$, we have for $i=1,2$,
  \begin{eqnarray*}
&&     \|I[N[u]]\|_{E_i \cap LE_i}\\&\le& C \sum_{|\al|= i-1}\|r^\de\<r\>^{ 1/2-\de'}\pa_x^\al N[u]\|_{L^2([0,\infty); L^2_x)}\\
       &\le& C (\|u\|_{E_1}^{1-s_1}\|u\|_{E_2}^{s_1}+\|u\|_{E_1}^{1-s_2}\|u\|_{E_2}^{s_2})^{p-1}
       \sum_{|\al|= i-1}\|r^{-\de}\<r\>^{-1/2+\de'}\pa_x^\al \pa u\|_{L^2([0,\infty); L^2_x)}\\
       &\le& C (\|u\|_{E_1}^{1-s_1}\|u\|_{E_2}^{s_1}+\|u\|_{E_1}^{1-s_2}\|u\|_{E_2}^{s_2})^{p-1} \|u\|_{LE_i}\ .
  \end{eqnarray*} This proves \eqref{eq-sub-inhomo}. A similar argument with \eqref{eq-Nonlinearity-rel2} instead of \eqref{eq-Nonlinearity-rel1} will yield \eqref{eq-sub-converg}.\end{prf}

With these two Propositions \ref{thm-sub-homo} and \ref{thm-sub-inho} in hand, it will be easy to show Theorem \ref{thm-sub}. Setting
$$\La_i:=\|u_0\|_{\dot H^i(\R^n)}+\|u_1\|_{\dot H^{i-1}(\R^n)}, \ i=1,2\ ,$$
we find by Propositions \ref{thm-sub-homo} and \ref{thm-sub-inho} that the mapping $\Phi$, defined by \eqref{eq-Picard}, is a contraction mapping from $X(2 C_1\La_1, 2C_1 \La_2; T)$ into itself, for any $T>0$ provided that
\beeq\label{eq-sub-small3}C_2  (2 C_1)^{p-1}(
\La_1^{1-s_1}\La_2^{s_1}+\La_1^{1-s_2}\La_2^{s_2})^{p-1}\le 1/2\ ,\eneq
and
\beeq\label{eq-sub-small2}C_3  (4 C_1)^{p-1}(\La_1^{1-s_1}\La_2^{s_1}+\La_1^{1-s_2}\La_2^{s_2})^{p-1}\le 1/2
\ .\eneq
Define a positive constant $C_0$ by
$$C_0^{-(p-1)}=\max(2 C_3  (4 C_1)^{p-1}, 2 C_2  (2 C_1)^{p-1})\ .$$
Then we see that when
\beeq\label{eq-sub-small}
\La_1^{1-s_1}\La_2^{s_1}+\La_1^{1-s_2}\La_2^{s_2}\le C_0\ ,
\eneq
 the map $\Phi$ is a contraction mapping of $X(2 C_1\La_1, 2C_1 \La_2; T)$ for any $T>0$, the global in time unique fixed point $u\in X(2 C_1\La_1, 2C_1 \La_2; \infty)$  is the solution which we seek.

To complete the proof of Theroem \ref{thm-sub}, we also need to establish the regularity of $u$, i.e., \beeq\label{eq-super-regularity}
\pt^i u\in C([0,\infty); H^{2-i}(\R^n)), i=0,1,
\eneq and the uniqueness of the solution $u$.

First, for the problem of regularity, it suffices to show
$$\pt^i u\in C([0,\infty); \dot H^{2-i}(\R^n)), i=0,1\ .$$
In fact, using the inequalities  \eqref{eq-sub-inhomo}, \eqref{eq-sub-small3} and the fact that $u\in LE_2$, we have
    \begin{eqnarray*}
     && \|\pa \pa_x (u(T)-u(0))\|_{L^2_x}\\
          &=& \|\pa \pa_x (\Phi[u](T)-u(0))\|_{L^2_x}\\
&\le&\|\pa \pa_x I[N[u]](T)\|_{L^2_x}+\|\pa \pa_x (u^{(0)}(T)-u(0))\|_{L^2_x}\nonumber\\
&\le &
 C_2 (\|u\|_{E_1}^{1-s_1}\|u\|_{E_2}^{s_1}+\|u\|_{E_1}^{1-s_2}\|u\|_{E_2}^{s_2})^{p-1} \\
 &&\times \|r^{-\de}\<r\>^{- 1/2+\de'} \pa \pa_x u\|_{L^2([0,T];L^2_x)}+o(1)\\
     &\le&  \|r^{-\de}\<r\>^{- 1/2+\de'} \pa \pa_x u\|_{L^2([0,T];L^2_x)}+o(1)=o(1) \ \nonumber
  \end{eqnarray*} as $T\rightarrow 0+$. This proves the continuity at $t=0$.
Recall that our solution satisfies $u=\Phi[u]$, which tells us that we can also view $u$ as the solution to the linear wave equation $\Box v = N[u](t_0 + t)$ with initial data $(u(t_0),\partial_t u(t_0))$ at any other time $t_0\in(0,\infty)$. Then a similar argument will give us the continuity at any $t\in [0,\infty)$.

Now, we turn to the proof of uniqueness. Assume there exists another solution $v\in X_{\infty}\cap C_t H^2\cap C_t^1 H^1$, with the same initial data.
Recall that $u,v\in C_t  H^2\cap C^1_t H^1$. If we restrict these solutions to small enough time interval $[0, T]$, owing to $\pa\pa_x^\al (u-v)(0)=0$, we have
\begin{eqnarray*}
\sum_{|\al|=i-1}  \|\pa \pa_x^\al v\|_{C([0,T]; L^2_x)}&\le &
\sum_{|\al|=i-1}  ( \|\pa\pa_x^\al (u-v)\|_{C([0,T]; L^2_x)}+\|\pa\pa_x^\al  u\|_{C([0,T]; L^2_x)})   \\
   &\le & o(1)+2 C_1 \La_i\ ,
\end{eqnarray*} as $T\rightarrow 0+$.
Using the inequality \eqref{eq-sub-converg}, we see that
\begin{eqnarray*}
&&  \|r^{-\de}\<r\>^{- 1/2+\de'} \pa (u-v)\|_{L^2([0,T];L^2_x)}\\
&=&\|r^{-\de}\<r\>^{- 1/2+\de'} \pa (\Phi[u]-\Phi[v])\|_{L^2([0,T];L^2_x)}\\
  & \le& C_3 \|r^{-\de}\<r\>^{- 1/2+\de'} \pa (u-v)\|_{L^2([0,T];L^2_x)}\\
  &&\times(
\|u\|_{E_1}^{1-s_1}\|u\|_{E_2}^{s_1}+\|u\|_{E_1}^{1-s_2}\|u\|_{E_2}^{s_2} +\|v\|_{E_1}^{1-s_1}\|v\|_{E_2}^{s_1}+\|v\|_{E_1}^{1-s_2}\|v\|_{E_2}^{s_2})^{p-1}\\
&\le & \frac 34   \|r^{-\de}\<r\>^{- 1/2+\de'} \pa (u-v)\|_{L^2([0,T];L^2_x)},
\end{eqnarray*} provided $T>0$ is small enough, where we have used \eqref{eq-sub-small} and \eqref{eq-sub-small2}.
By this we arrive at the conclusion that $u=v$ for $t\in [0,T]$, which shows the uniqueness. This completes the proof of Theorem \ref{thm-sub}.

\subsection{Glassey conjecture when $p=p_c$ and $n\ge 3$}
Consider \eqref{eq-NLW} for $p=p_c$ and $n\ge 3$.

The estimate of the homogeneous solution, $u^{(0)}$, is given by Proposition \ref{thm-sub-homo}.
We only need to give the estimate for the inhomogeneous part.
\begin{prop}\label{thm-crit-inho} Let $p=p_c$ and $n\ge 3$, $u\in X_T$ and $s\in (1/2,1]$. Define
 \beeq\label{eq-crit-delta}\de=\frac{n-2 s}4 (p-1)\ .\eneq
 Then there is a positive constant $C_4$, independent of $T>0$, such that the following estimates hold
\beeq\|I[N[u]]\|_{E_i\cap LE_i}\le C_4 \log(2+T) (\|u\|_{E_1}^{1/2}\|u\|_{E_2}^{1/2}+\|u\|_{E_1}^{1-s}\|u\|_{E_2}^{s})^{p-1} \|u\|_{LE_i},\ i=1,2\ .\label{eq-crit-inhomo}
\eneq
  Moreover, we have
\begin{eqnarray}
&&\|\Phi[u]-\Phi[v]\|_{E_1\cap LE_1} \label{eq-crit-converg}\\
&  &\le C_5 \log(2+T) \|u-v\|_{LE_1} \nonumber\\
&&\times (\|u\|_{E_1}^{1/2}\|u\|_{E_2}^{1/2}+\|u\|_{E_1}^{1-s}\|u\|_{E_2}^{s} +\|v\|_{E_1}^{1/2}\|v\|_{E_2}^{1/2}+\|v\|_{E_1}^{1-s}\|v\|_{E_2}^{s})^{p-1}\ ,  \nonumber
\end{eqnarray}
 for some $C_5$.
\end{prop}

\begin{prf}
First, by \eqref{eq-sub-decay} with $s_1=1/2$ and $s_2=s$, we have
\beeq\label{eq-crit-decay} \|u(r\omega)\|_{L^2_\omega}\le C r^{s-n/2}\<r\>^{1/2-s}(\|u\|_{L^2_x}^{1/2}\|\pa_r u\|_{L^2_x}^{1/2}+\|u\|_{L^2_x}^{1-s}\|\pa_r u\|_{L^2_x}^{s})\ .\eneq
By Lemma \ref{thm-radial}, we have for $u\in X_T$,
  \beeq\label{eq-crit-decay2}
    |\pa u|\le C r^{s-n/2}\<r\>^{1/2-s}(\|u\|_{E_1}^{1/2}\|u\|_{E_2}^{1/2}+\|u\|_{E_1}^{1-s}\|u\|_{E_2}^{s})\ .
  \eneq
From \eqref{eq-crit-decay2} and \eqref{eq-Nonlinearity-rel1}, it is clear that, for $i=1,2$,
  \begin{eqnarray*}
 &&   \sum_{|\al|= i-1}\|r^\de\<r\>^{ 1/2-\de}\pa_x^\al N[u]\|_{L^2_x}\\
  &\le& C \sum_{|\al|= i-1}\|r^\de\<r\>^{1/2-\de}|\pa u|^{p-1}\pa_x^\al \pa u\|_{L^2_x} \\
          &\le & C  (\|u\|_{E_1}^{1/2}\|u\|_{E_2}^{1/2}+\|u\|_{E_1}^{1-s}\|u\|_{E_2}^{s})^{p-1}\\
          &&\times \sum_{|\al|= i-1}\|r^{\de+(s-n/2)(p-1)}\<r\>^{ 1/2-\de+(1/2-s)(p-1)}\pa_x^\al \pa u\|_{L^2_x}\\
          &= & C (\|u\|_{E_1}^{1/2}\|u\|_{E_2}^{1/2}+\|u\|_{E_1}^{1-s}\|u\|_{E_2}^{s})^{p-1}\sum_{|\al|= i-1}\|r^{-\de}\<r\>^{-1/2+\de}\pa_x^\al \pa u\|_{L^2_x}\ .\\
  \end{eqnarray*}
Now applying Lemma \ref{thm-KSS-HY-inh} to $\pa_x^\al u$ with $|\al|\le 1$ and $u_0=u_1=0$, we have for $i=1,2$,
  \begin{eqnarray*}
&&     \|I[N[u]]\|_{E_i \cap LE_i}\\&\le& C(\log(2+T))^{1/2} \sum_{|\al|= i-1}\|r^\de\<r\>^{1/2-\de}\pa_x^\al N[u]\|_{L^2([0,T]; L^2_x)}\\
       &\le& C(\log(2+T))^{1/2} (\|u\|_{E_1}^{1/2}\|u\|_{E_2}^{1/2}+\|u\|_{E_1}^{1-s}\|u\|_{E_2}^{s})^{p-1}\\
&&\times       \sum_{|\al|= i-1}\|r^{-\de}\<r\>^{-1/2+\de}\pa_x^\al \pa u\|_{L^2([0,T]; L^2_x)}\\
       &\le& C \log(2+T) (\|u\|_{E_1}^{1/2}\|u\|_{E_2}^{1/2}+\|u\|_{E_1}^{1-s}\|u\|_{E_2}^{s})^{p-1} \|u\|_{LE_i}\ .
  \end{eqnarray*} This proves \eqref{eq-crit-inhomo}. A similar argument with \eqref{eq-Nonlinearity-rel2} instead of \eqref{eq-Nonlinearity-rel1} will yield \eqref{eq-crit-converg}.
  Here, for later use, we record the following inequality which is a direct consequence of the last one,
\begin{eqnarray}
&&  \sup_{t\in [0,T]} \|\pa \pa_x I[N[u]](t)\|_{L^2_x(\R^n)}\label{eq-crit-regularity}\\
&\le& C_4(\log(2+T))^{1/2}  \|r^{-\de}\<r\>^{-1/2+\de}\pa_x \pa u\|_{L^2([0,T]; L^2_x)} \nonumber\\
  &&\times  (\|u\|_{E_1}^{1/2}\|u\|_{E_2}^{1/2}+\|u\|_{E_1}^{1-s}\|u\|_{E_2}^{s})^{p-1}\ .\nonumber
\end{eqnarray}
  \end{prf}

With these two Propositions \ref{thm-sub-homo} and \ref{thm-crit-inho} in hand, the proof of Theorem \ref{thm-crit} proceeds similarly to that of Theorem \ref{thm-sub}. With $$\La_i:=\|u_0\|_{\dot H^i(\R^n)}+\|u_1\|_{\dot H^{i-1}(\R^n)}, \ i=1,2\ ,$$
we find by Propositions \ref{thm-sub-homo} and \ref{thm-crit-inho} that the mapping $\Phi$, defined by \eqref{eq-Picard}, is a contraction mapping from $X(2 C_1\La_1, 2C_1 \La_2; T)$ into itself, for any $T>0$ provided that
\beeq\label{eq-crit-small3}C_4 \log (2+T)  (2 C_1)^{p-1} (
\La_1^{1/2}\La_2^{1/2}+\La_1^{1-s}\La_2^{s})^{p-1}\le 1/2\ ,\eneq
and
\beeq\label{eq-crit-small2}C_5 \log (2+T)  (4 C_1)^{p-1} (
\La_1^{1/2}\La_2^{1/2}+\La_1^{1-s}\La_2^{s})^{p-1}\le 1/2
\ .\eneq
Define a positive constant $C_6$ by
$$C_6^{-1}=\max(2 C_5 (4 C_1)^{p-1}, 2 C_4  (2 C_1)^{p-1})\ ,$$
and set $T_*$ according to
$$\log (2+T_*)   (
\La_1^{1/2}\La_2^{1/2}+\La_1^{1-s}\La_2^{s})^{p-1}=C_6\ ,$$
which is possible in general only if $$\ep=\La_1^{1/2}\La_2^{1/2}+\La_1^{1-s}\La_2^{s}\ll 1\ .$$
That is
\beeq\label{eq-crit-lifespan}
 T_*=\exp(C_6 \ep^{1-p})-2,\ \ep\ll 1\ .
\eneq
Since $\Phi$ is a  contraction mapping in $X(2 C_1\La_1, 2C_1 \La_2; T_*)$, the unique fixed point $u\in X(2 C_1\La_1, 2C_1 \La_2; T_*)$  is the solution which we seek.

To complete the proof of Theroem \ref{thm-crit}, we need also to establish the regularity of $u$, i.e., \beeq\label{eq-super-regularity}
\pt^i u\in C([0,T_*]; H^{2-i}(\R^n)), i=0,1,
\eneq and the uniqueness of the solution $u$.

First, for the problem of regularity, it suffices to show
$$\pt^i u\in C([0,T_*]; \dot H^{2-i}(\R^n)), i=0,1\ .$$
Indeed, since $u\in LE_2(T_*)$, we know that $$\|r^{-\de}\<r\>^{-1/2+\de}\pa_x \pa u\|_{L^2([0,T_*]; L^2_x)}<\infty\ ,$$
and so
$$\lim_{T\rightarrow 0+}\|r^{-\de}\<r\>^{- 1/2+\de}\pa_x \pa u\|_{L^2([0,T]; L^2_x)}=0\ .$$
 Using the inequality  \eqref{eq-crit-regularity} and \eqref{eq-crit-small3}, we have
    \begin{eqnarray*}
     && \|\pa \pa_x (u(T)-u(0))\|_{L^2_x}\\
          &=& \|\pa \pa_x (\Phi[u](T)-u(0))\|_{L^2_x}\\
&\le&\|\pa \pa_x I[N[u]](T)\|_{L^2_x}+\|\pa \pa_x (u^{(0)}(T)-u(0))\|_{L^2_x}\nonumber\\
     &\le&
     C_4(\log(2+T))^{1/2}  \|r^{-\de}\<r\>^{-1/2+\de}\pa_x \pa u\|_{L^2([0,T]; L^2_x)} \nonumber\\
  &&\times  (\|u\|_{E_1}^{1/2}\|u\|_{E_2}^{1/2}+\|u\|_{E_1}^{1-s}\|u\|_{E_2}^{s})^{p-1}\nonumber+o(1)\\
  &\le & \|r^{-\de}\<r\>^{-1/2+\de}\pa_x \pa u\|_{L^2([0,T]; L^2_x)}+o(1)  =o(1) \ \nonumber
  \end{eqnarray*} as $T\rightarrow 0+$. This proves the continuity at $t=0$. A similar argument will give us the continuity at any $t\in [0,T_*]$.

Now, we turn to the proof of uniqueness. Assume that there exists another solution $v\in X_{T_*}\cap C_t H^2\cap C_t^1 H^1$, with the same initial data.
Recall that $u,v\in C_t  H^2\cap C^1_t H^1$. If we restrict these solutions to small enough time interval $[0, T]$, owing to $\pa\pa_x^\al (u-v)(0)=0$, we have
\begin{eqnarray*}
\sum_{|\al|=i-1}  \|\pa \pa_x^\al v\|_{C([0,T]; L^2_x)}&\le &
\sum_{|\al|=i-1}  ( \|\pa\pa_x^\al (u-v)\|_{C([0,T]; L^2_x)}+\|\pa\pa_x^\al  u\|_{C([0,T]; L^2_x)})   \\
   &\le & o(1)+2 C_1 \La_i\ ,
\end{eqnarray*} as $T\rightarrow 0+$.
Using the inequality \eqref{eq-sub-converg}, we see that
\begin{eqnarray*}
&&  \|r^{-\de}\<r\>^{-1/2+\de} \pa (u-v)\|_{L^2([0,T];L^2_x)}\\
&=&\|r^{-\de}\<r\>^{-1/2+\de} \pa (\Phi[u]-\Phi[v])\|_{L^2([0,T];L^2_x)}\\
  & \le& C_5 \log(2+T) \|r^{-\de}\<r\>^{-1/2+\de} \pa (u-v)\|_{L^2([0,T];L^2_x)}\\
  &&\times(
\|u\|_{E_1}^{1/2}\|u\|_{E_2}^{1/2}+\|u\|_{E_1}^{1-s}\|u\|_{E_2}^{s} +\|v\|_{E_1}^{1/2}\|v\|_{E_2}^{1/2}+\|v\|_{E_1}^{1-s}\|v\|_{E_2}^{s})^{p-1}\\
&\le & \frac 34   \|r^{-\de}\<r\>^{-1/2+\de} \pa (u-v)\|_{L^2([0,T];L^2_x)},
\end{eqnarray*} provided $T>0$ is small enough, where we have used \eqref{eq-crit-small2}.
By this, we conclude that $u=v$ for $t\in [0,T]$, which shows the uniqueness. This completes the proof of Theorem \ref{thm-crit}.

\section{Glassey conjecture when $p<p_c$ and $n\ge 2$}

 In this section, we aim at giving the proof of Theorem \ref{thm-super} for $n\ge 2$. As we will see, the argument in the previous section can be adapted to the scale-supercritical case $p<p_c$, for $n\ge 3$. The argument in the previous section does not apply when $n=2$, owing to the fact that current techniques do not yield the inhomogeneous KSS type estimates \eqref{eq-KSS-HY-inh} for $n=2$. Alternatively, applying the homogeneous estimates in Lemma \ref{thm-KSS-HY} gives us the proof.

In this section, by $\de$ in $LE$ norm, we mean
\beeq\label{eq-super-delta}\de=\left\{\begin{array}{ll}
  \frac{(n-1)(p-1)}2,& 1<p<1+\frac{1}{n-1}\\
  \frac{(n-1)(p-1)}4,& 1+\frac{1}{n-1}\le p<1+\frac{2}{n-1}=p_c\ .
\end{array}\right.\eneq Note that $0<\de<1/2$.

We aim at showing that $\Phi$ is a contraction mapping of $X(R_1, R_2; T)$, if we choose $R_1$, $R_2$ and $T$ suitably. As before, the estimate of the homogeneous solution, $u^{(0)}$, is given by Proposition \ref{thm-sub-homo}.
We only need to obtain a similar estimate for the inhomogeneous part.
\begin{prop}\label{thm-super-inho} Let $1<p<p_c$, $u\in X_T$ and $\de$ as in \eqref{eq-super-delta}. Then
  there is a positive constant $C_7$, independent of $T>0$, such that the following estimates hold
\begin{eqnarray}
  &&   \|I[N[u]]\|_{E_i}+\|I[N[u]]\|_{LE_i}\label{eq-super-inhomo}\\
&& \le  C_7 T^{1-(n-1)(p-1)/2} (\|u\|_{E_1\cap LE_1}\|u\|_{E_2\cap LE_2})^{(p-1)/2}\|u\|_{LE_i},\ i=1,2\ . \nonumber
\end{eqnarray}
  Moreover, we have
\begin{eqnarray}
&&\|\Phi[u]-\Phi[v]\|_{E_1\cap LE_1} \label{eq-super-converg}\\
&\le & C_8 T^{1-(n-1)(p-1)/2}(\|u\|_{E_1\cap LE_1}\|u\|_{E_2\cap LE_2}+\|v\|_{E_1\cap LE_1}\|v\|_{E_2\cap LE_2})^{(p-1)/2}  \nonumber\\
&&\times \|u-v\|_{LE_1}\nonumber
\end{eqnarray}
 for some $C_8$, independent of $T>0$.
\end{prop}
\begin{rem}\label{thm-super-uniq}
  From the proof of \eqref{eq-super-converg}, we can extract the following estimates. If $1<p<\min(p_c,2)$, then
  \begin{eqnarray}
     && \|r^{-\de}\pa(\Phi[u]-\Phi[v])\|_{L^2([0,T];L^2_x)} \label{eq-super-uniq1}\\
     &\le& C_8 T^{1-(n-1)(p-1)/2} (\|u\|_{E_1}\|u\|_{E_2}+\|v\|_{E_1}\|v\|_{E_2})^{(p-1)/2} \nonumber\\
     &&\times      \|r^{-\de}\pa(u-v)\|_{L^2([0,T];L^2_x)}\ .\nonumber
  \end{eqnarray}
If $2\le p<3$ and $n=2$,
\begin{eqnarray}
  &&\|r^{-\de}(\Phi[u]-\Phi[v])\|_{L^2([0,T];L^2_x)}\label{eq-super-uniq2}\\
  &\le&
  C_8 T^{(3-p)/4} (\|u\|_{E_1}\|u\|_{E_2}+\|v\|_{E_1}\|v\|_{E_2})^{(p-2)/2}\nonumber\\
&&\times\left(\|r^{-\de}\pa u\|_{L^2([0,T];L^2_x)}\|r^{-\de}\pa \pa_x u\|_{L^2([0,T];L^2_x)}\right.\nonumber\\
&&+\left.\|r^{-\de}\pa v\|_{L^2([0,T];L^2_x)}\|r^{-\de}\pa \pa_x v\|_{L^2([0,T];L^2_x)}\right)^{1/2}\nonumber\\
&&\times\|r^{-\de}\pa (u-v)\|_{L^2([0,T];L^2_x)}\ .\nonumber
\end{eqnarray}
\end{rem}
\begin{rem}\label{thm-super-regu}
  From the proof of \eqref{eq-super-inhomo} with $i=2$, we can extract the following estimates. If   $1<p<\min(p_c,2)$, then
  \begin{eqnarray}
     && \|\pa \pa_x (\Phi[u](T)-u(0))\|_{L^2_x} \label{eq-super-regu1}\\
&\le&\|\pa \pa_x I[N[u]](T)\|_{L^2_x}+\|\pa \pa_x (u^{(0)}(T)-u(0))\|_{L^2_x}\nonumber\\
     &\le& C_7 T^{1-(n-1)(p-1)/2} (\|u\|_{E_1}\|u\|_{E_2})^{(p-1)/2} T^{\de-1/2} \|r^{-\de}\pa \pa_x u\|_{L^2([0,T];L^2_x)}+o(1) \ \nonumber
  \end{eqnarray} as $T\rightarrow 0+$.
If $2\le p<3$ and $n=2$, \begin{eqnarray}
  &&\|\pa \pa_x (\Phi[u](T)-u(0))\|_{L^2_x}\label{eq-super-regu2}\\
  &\le&\|\pa \pa_x I[N[u]](T)\|_{L^2_x}+\|\pa \pa_x (u^{(0)}(T)-u(0))\|_{L^2_x}\nonumber\\
  &\le&
  C_7   (\|u\|_{E_1}\|u\|_{E_2})^{(p-2)/2}  \|r^{-(p-1)/4}\pa u\|_{L^2([0,T]; L^2_x)}^{1/2}\nonumber\\
  &&\times
  \|r^{-(p-1)/4}\pa \pa_x u\|_{L^2([0,T]; L^2_x)}^{3/2}+o(1)\nonumber
   \end{eqnarray}
   as $T\rightarrow 0+$.
\end{rem}

\begin{prf}
  We will deal with three different cases: $1<p<1+1/(n-1)$ when $n\ge 2$; $2\le p<3$ when $n=2$; and $1+1/(n-1)\le p<p_c$ when $n\ge 3$.

  {\bf Case i) $1<p<1+1/(n-1)$ with $\de=(n-1)(p-1)/2$}. First, by \eqref{eq-trace-edpt} and Lemma \ref{thm-radial}, we have for $u\in X_T$,
  \beeq\label{eq-decay}
    |\pa u|\le C r^{-(n-1)/2}(\|u\|_{E_1}\|u\|_{E_2})^{1/2}\ .
  \eneq
  By \eqref{eq-decay} and \eqref{eq-Nonlinearity-rel1}, it is clear that, for $i=1,2$,
  \begin{eqnarray*}
    \sum_{|\al|= i-1}\|\pa_x^\al N[u]\|_{L^2_x} &\le& C \sum_{|\al|= i-1}\||\pa u|^{p-1}\pa_x^\al \pa u\|_{L^2_x} \\
          &\le & C  (\|u\|_{E_1}\|u\|_{E_2})^{(p-1)/2}\sum_{|\al|= i-1}\|r^{-(n-1)(p-1)/2}\pa_x^\al \pa u\|_{L^2_x}\\
          &= & C  (\|u\|_{E_1}\|u\|_{E_2})^{(p-1)/2}\sum_{|\al|= i-1}\|r^{-\de}\pa_x^\al \pa u\|_{L^2_x}\ .\\
  \end{eqnarray*}
  Now applying Lemma \ref{thm-KSS-HY} to $\pa_x^\al u$ with $|\al|\le 1$ and $u_0=u_1=0$, we have for $i=1,2$,
  \begin{eqnarray*}
     &&\|I[N[u]]\|_{E_i \cap LE_i}\\
     &\le& C \sum_{|\al|= i-1}\|\pa_x^\al N[u]\|_{L^1([0,T]; L^2_x)}\\
       &\le& C T^{1/2} \sum_{|\al|= i-1}\|\pa_x^\al N[u]\|_{L^2([0,T]; L^2_x)}\\
       &\le& C T^{1/2}  (\|u\|_{E_1}\|u\|_{E_2})^{(p-1)/2} \sum_{|\al|= i-1}\|r^{-\de}\pa_x^\al \pa u\|_{L^2([0,T]; L^2_x)}\\
       &\le& C T^{1-\de}  (\|u\|_{E_1}\|u\|_{E_2})^{(p-1)/2} T^{\de-1/2} \sum_{|\al|= i-1}\|r^{-\de}\pa_x^\al \pa u\|_{L^2([0,T]; L^2_x)}\\
       &\le &C T^{1-\de}  (\|u\|_{E_1}\|u\|_{E_2})^{(p-1)/2}\|u\|_{LE_i}\ .
  \end{eqnarray*} This proves \eqref{eq-super-inhomo}. A similar argument with \eqref{eq-Nonlinearity-rel2} instead of \eqref{eq-Nonlinearity-rel1} will yield \eqref{eq-super-converg}.

  {\bf Case ii) $1+1/(n-1)\le p<p_c$, $n\ge 3$ and $\de=(n-1)(p-1)/4$}. In this case, we may use Lemma \ref{thm-KSS-HY-inh} instead. Applying it to $\pa_x^\al u$ with $|\al|\le 1$ and $u_0=u_1=0$, we have for $i=1,2$,
  \begin{eqnarray*}
  &&\|I[N[u]]\|_{E_i}+\|I[N[u]]\|_{LE_i}\\
     &\le& C T^{1/2-\de} \sum_{|\al|= i-1}\|r^\de \pa_x^\al N[u]\|_{L^2([0,T]; L^2_x)}\\
     &\le& CT^{1/2-\de}  (\|u\|_{E_1}\|u\|_{E_2})^{(p-1)/2} \sum_{|\al|= i-1}\|r^{-(n-1)(p-1)/2} r^{\de} \pa_x^\al \pa u\|_{L^2([0,T]; L^2_x)}\\
       &\le& C T^{1-2\de}  (\|u\|_{E_1}\|u\|_{E_2})^{(p-1)/2} T^{\de-1/2} \sum_{|\al|= i-1}\|r^{-\de}\pa_x^\al \pa u\|_{L^2([0,T]; L^2_x)}\\
       &\le &C T^{1-2\de}  (\|u\|_{E_1}\|u\|_{E_2})^{(p-1)/2}\|u\|_{LE_i}\ ,
  \end{eqnarray*} where we have used \eqref{eq-decay} and \eqref{eq-Nonlinearity-rel1}.
Using \eqref{eq-Nonlinearity-rel2} instead of \eqref{eq-Nonlinearity-rel1}, \eqref{eq-super-converg} follows similarly.

  {\bf Case iii) $2\le p<3$,  $n=2$ and $\de=(p-1)/4$}.
Notice that Lemma \ref{thm-trace-variant} with $s=(3-p)/4>0$ gives us
$$\|r^{-(p-3)/4} u\|_{L^\infty_r L^2_\omega}\le \sqrt{2} \|r^{-(p-1)/4}u\|_{L^2_x(\R^2)}^{1/2}\|r^{-(p-1)/4}\pa_x u\|_{L^2_x(\R^2)}^{1/2}\ .
$$
Then for $i=1,2$, we obtain by using Lemma \ref{thm-KSS-HY}
  \begin{eqnarray*}
      && \|I[N[u]]\|_{E_i}+\|I[N[u]]\|_{LE_i}\\
      &\le& C \sum_{|\al|= i-1}\| \pa_x^\al N[u]\|_{L^1([0,T]; L^2_x(\R^2))}\\
      &\le& C  (\|u\|_{E_1}\|u\|_{E_2})^{(p-2)/2} \sum_{|\al|= i-1}\| r^{-(p-2)/2} \pa u \pa_x^\al \pa u\|_{L^1([0,T];L^2_x(\R^2))}\\
      &\le& C  (\|u\|_{E_1}\|u\|_{E_2})^{(p-2)/2}  \|r^{-(p-3)/4}\pa u\|_{L^2([0,T]; L^\infty_x)} \\
      &&\times \sum_{|\al|= i-1}\|r^{-(p-1)/4}\pa_x^\al \pa u\|_{L^2([0,T];L^2_x)}\\
            &\le& C  (\|u\|_{E_1}\|u\|_{E_2})^{(p-2)/2}  \|r^{-(p-3)/4}\pa u\|_{L^2([0,T]; L^\infty_r L^2_\omega)} \\
      &&\times \sum_{|\al|= i-1}\|r^{-(p-1)/4}\pa_x^\al \pa u\|_{L^2([0,T];L^2_x)}\\
      &\le &C T^{1-2\de}  (\|u\|_{E_1}\|u\|_{E_2})^{(p-2)/2}(\|u\|_{LE_1}\|u\|_{LE_2})^{1/2} \|u\|_{LE_i}\ ,
  \end{eqnarray*} where we have used \eqref{eq-Nonlinearity-rel1} and Lemma \ref{thm-radial}.

The estimate \eqref{eq-super-converg} follows from the similar arguments by using \eqref{eq-Nonlinearity-rel2} instead of \eqref{eq-Nonlinearity-rel1}. This completes the proof. \end{prf}

With these two Propositions \ref{thm-sub-homo} and \ref{thm-super-inho} in hand, it will be easy to show Theorem \ref{thm-super}. Setting
$$\La_i:=\|u_0\|_{\dot H^i(\R^n)}+\|u_1\|_{\dot H^{i-1}(\R^n)}, \ i=1,2\ ,$$
we find by Propositions \ref{thm-sub-homo} and \ref{thm-super-inho} that the mapping $\Phi$, defined by \eqref{eq-Picard}, is a contraction mapping from $X(2 C_1\La_1, 2C_1 \La_2; T)$ into itself, provided that
$$C_7 T^{1-(n-1)(p-1)/2} (2 C_1)^{p-1}(\La_1\La_2)^{(p-1)/2}\le 1/2$$
and
$$C_8 T^{1-(n-1)(p-1)/2} (4 C_1)^{p-1}(\La_1\La_2)^{(p-1)/2}\le 1/2\ .$$
Define a positive constant $C_9$ by
$$C_9^{-\left(1-(n-1)(p-1)/2\right)}=\max(2 C_8  (4 C_1)^{p-1}, 2 C_7  (2 C_1)^{p-1})\ ,$$
and set $T_*$ according to
$$C_9^{-\left(1-(n-1)(p-1)/2\right)} T_*^{1-(n-1)(p-1)/2} (\La_1\La_2)^{(p-1)/2}=1\ ,$$
that is
\beeq\label{eq-super-lifespan}
 T_*=C_9 (\La_1\La_2)^{-\frac{p-1}{2-(n-1)(p-1)}}\ .
\eneq
Since $\Phi$ is a contraction mapping of $X(2 C_1\La_1, 2C_1 \La_2; T_*)$, the unique fixed point $u\in X(2 C_1\La_1, 2C_1 \La_2; T_*)$  is the solution which we seek.

To complete the proof of Theroem \ref{thm-super}, we also need to establish the uniqueness of $u$ in $X_{T_*}$, and the regularity of $u$, i.e., \beeq\label{eq-super-regularity}
\pt^i u\in C([0,T_*]; H^{2-i}(\R^n)), i=0,1.
\eneq
First, for the proof of uniqueness, assume that there exists another solution $v\in X_{T_*}$, with the same initial data.
Recall the estimates \eqref{eq-super-uniq1} and \eqref{eq-super-uniq2}. If we restrict these solutions to small enough $0<T<T_*$, we have
\begin{eqnarray*}
  \|r^{-\de}\pa (u-v)\|_{L^2([0,T]; L^2_x)}&=&
   \|r^{-\de}\pa (\Phi[u]-\Phi[v])\|_{L^2([0,T]; L^2_x)}\\
   &\le &\frac 12 \|r^{-\de}\pa (u-v)\|_{L^2([0,T]; L^2_x)}\ .
\end{eqnarray*} Combining this with the fact that $u$ and $v$ share the same initial data, we conclude that $u=v$ for $t\in [0,T]$, which shows the uniqueness.

By Remark \ref{thm-super-regu}, $1-(n-1)(p-1)/2+\delta-1/2>0$, and the fact that $r^{-\de}\pa\pa_x u\in L^2([0,T_*]; L^2_x)$ (and so $$\|r^{-\de}\pa\pa_x u\|_{L^2([0,T]; L^2_x)}=o(1)$$ as $T\rightarrow 0+$), we see that $\pa\pa_x u(T)$ converges to $\pa\pa_x u(0)$ in $L^2_x$. This tells us that the continuity at $t=0$. A similar argument will give us the continuity at any $t\in [0,T_*]$. This completes the proof of Theorem \ref{thm-super}.

\section{Glassey conjecture when $n=2$, $p> p_c$}

For this case, it seems not enough for us to prove global results by applying the KSS type estimates, mainly because we do not have the favorable inhomogeneous estimates as \eqref{eq-KSS-HY-inh} in Lemma \ref{thm-KSS-HY-inh}.

Instead, we want to present a proof based on the recent generalized Strichartz estimates of Smith, Sogge and Wang \cite{SSW10} (with the previous radial estimates in Fang and Wang \cite{FW06}).
\begin{lem}[Generalized Strichartz estimates]\label{thm-SSW}
Let $n=2$ and $q\in (2,\infty)$.  For any solution $u=u(t,x)$ to the wave equation \eqref{eq-lin-wave}, we have the following inequality with $s=1-1/q$,
\beeq\label{eq-SSW}
    \|\pa u\|_{L^q ([0,\infty); L^\infty_r L^2_\omega (\R^2))}\le C_q (\|\pa_x u_0\|_{\dot H^{s}_x}+\|u_1\|_{\dot H^s_x}+\|F\|_{L^1_t \dot H^{s}_x})\ ,
\eneq where $C$ is independent of the functions $u_0$, $u_1$ and $F$.
\end{lem}

With these estimates, we are able to present a simple proof of Theorem \ref{thm-sub-n2}.
Let
$$\La_i:=\|u_0\|_{\dot H^i(\R^n)}+\|u_1\|_{\dot H^{i-1}(\R^n)}, \ i=1,2\ .$$
By using \eqref{eq-Nonlinearity-rel1} and the energy estimates, we have
\begin{eqnarray}
  \|\pa \Phi[u]\|_{L^\infty_t L^2_x}&\le & C\Lambda_1+C\|N[u]\|_{L^1_t L^2_x}\label{eq-2D-BD}\\
     &\le & C\Lambda_1+C\|\pa u\|_{L^{p-1}_t L^\infty_x}^{p-1} \|\pa u\|_{L^\infty_t L^2_x} \ ,\nonumber
\end{eqnarray}and
\begin{eqnarray}
  \|\pa \pa_x \Phi[u]\|_{L^\infty_t L^2_x}&\le & C\Lambda_2+C\|\pa_x N[u]\|_{L^1_t L^2_x}\label{eq-2D-BD2}\\
     &\le & C\Lambda_2+C\|\pa u\|_{L^{p-1}_t L^\infty_x}^{p-1} \|\pa \pa_x u\|_{L^\infty_t L^2_x} \ .\nonumber
\end{eqnarray}
Recall the convex inequality
$$\|f\|_{\dot H^{1-\theta}}\le\|f\|_{L^2}^{\theta}\|f\|_{\dot H^1}^{1-\theta}, \ \theta\in [0,1],$$
together with \eqref{eq-Nonlinearity-rel1},  Lemma \ref{thm-SSW} and Lemma \ref{thm-radial}. We see that for $p>3$,
\begin{eqnarray}
&& \|\pa \Phi[u]\|_{L^{p-1}_t L^\infty_x}\label{eq-2D-BD3}\\
&\le & \|\pa \Phi[u]\|_{L^{p-1}_t L^\infty_r L^2_\omega}\nonumber\\
    &\le&   C\|\pa u(0)\|_{\dot H^{1-1/(p-1)}}+C\| N[u]\|_{L^1_t \dot H^{1-1/(p-1)}}\nonumber\\
     &\le & C\Lambda_1^{1/(p-1)}\La_2^{1-1/(p-1)}+C\|\pa u\|_{L^{p-1}_t L^\infty_x}^{p-1}
    \|\pa  u\|_{L^\infty_t L^2_x}^{1/(p-1)}
     \|\pa \pa_x u\|_{L^\infty_t L^2_x}^{1-1/(p-1)} \ .\nonumber
\end{eqnarray}
Moreover, we have
\begin{eqnarray}
&&\|\pa (\Phi[u]-\Phi[v])\|_{L^\infty_t L^2_x} \label{eq-2D-converg}\\
 &\le & C \|N[u]-N[v]\|_{L^1_t L^2_x}\nonumber\\
     &\le & C(\|\pa u\|_{L^{p-1}_t L^\infty_x}+\|\pa v\|_{L^{p-1}_t L^\infty_x})^{p-1} \|\pa (u-v)\|_{L^\infty_t L^2_x} \ .\nonumber
\end{eqnarray}

Let $\ep_0>0$ be the number such that $$C (4C\ep_0)^{p-1}= 1/2\ .$$
If
$$\ep= \Lambda_1^{1/(p-1)}\La_2^{1-1/(p-1)}\le \ep_0\ ,$$
then we see that $\Phi$ is a contraction mapping in $Y(2 C \La_1, 2 C \La_2, 2C \ep)$. Here the complete space $Y(R_1, R_2, R_3)$ is defined as
 \begin{eqnarray*}
&   Y(R_1, R_2, R_3)=&\{u\in C_t H_{\rm{rad}}^1\cap C^1_t L_{\rm{rad}}^2;  \|\pa u\|_{L^\infty_t L^2_x}\le R_1,\\
&& \|\pa \pa_x u\|_{L^\infty_t L^2_x}\le R_2,\
 \|\pa u\|_{L^{p-1}_t L^\infty_x}\le R_3
 \} \end{eqnarray*} with the metric $\rho(u,v)=\|\pa (u-v)\|_{L^\infty_t L^2_x}$.

To prove the regularity, we only need to show the continuity at $t=0$. For that, since $\pa u\in L^{p-1}([0,\infty); L^\infty_x)$, we have
\begin{eqnarray*}
&&\|\pa\pa_x (u(t)-u(0))\|_{ L^2_x}\\
 &\le & \|\pa\pa_x I[N[u]](t)\|_{ L^2_x}+\|\pa\pa_x (u^{(0)}(t)-u(0))\|_{ L^2_x}\nonumber\\
 &\le & \|\pa_x N[u]\|_{L^1([0,t]; L^2_x)}+o(1)\nonumber\\
     &\le & C\|\pa u\|_{L^{p-1}([0,t]; L^\infty_x)}^{p-1}\|\pa\pa_x u\|_{L^\infty_t L^2_x} +o(1)=o(1)\nonumber
\end{eqnarray*} as $t\rightarrow 0+$. This tells us that $u\in C_t \dot H^2 \cap C_t^1 \dot H^1$.

For uniqueness, suppose that there exists another solution $v\in Y\cap C_t H^2\cap C_t^1 H^1$, with the same initial data.
Using the inequality \eqref{eq-2D-converg}, we see that
\begin{eqnarray*}
&&\|\pa (u-v)\|_{C_t([0,T]; L^2_x)} \\
&=&\|\pa (\Phi[u]-\Phi[v])\|_{C_t([0,T]; L^2_x)} \\
     &\le & C(\|\pa u\|_{L^{p-1}_t([0,T]; L^\infty_x)}+\|\pa v\|_{L^{p-1}_t ([0,T]; L^\infty_x)})^{p-1} \|\pa (u-v)\|_{C_t([0,T]; L^2_x)} \\
     &\le & o(1) \|\pa (u-v)\|_{C_t([0,T]; L^2_x)}
\end{eqnarray*} as $T\rightarrow 0+$. Thus by choosing $T>0$ small enough, we  conclude that $u=v$ for $t\in [0,T]$, which shows the uniqueness. This completes the proof of Theorem \ref{thm-sub-n2}.

\end{document}